\documentclass[11pt]{amsart}

\usepackage{amsmath,amsthm}
\usepackage{amssymb}

\input xy
\xyoption{all}

\CompileMatrices

\DeclareMathOperator{\sln}{\mathfrak{sl}}

\DeclareMathOperator{\tr}{tr}
\DeclareMathOperator{\Tr}{Tr}
\DeclareMathOperator{\iden}{1}
\DeclareMathOperator{\im}{Im}

\DeclareMathOperator{\End}{End}
\DeclareMathOperator{\Hom}{Hom}
\DeclareMathOperator{\rep}{Rep}
\DeclareMathOperator{\spn}{span}
\DeclareMathOperator{\ab}{Ab}

\DeclareMathOperator{\add}{Add}

\DeclareMathOperator{\Ob}{Ob}

\DeclareMathOperator{\idem}{Idem}
\DeclareMathOperator{\mult}{mult}

\renewcommand{\1}{{{\rm 1 \hspace{-0.35em} 1}}}
\newcommand{\smone}{{{\rm 1 \hspace{-0.25em} 1}}}
\newcommand{\ident}{{{\rm 1 \hspace{-0.35em} 1}}}
\newcommand{\bmx}{\left(\begin{array}}
\newcommand{\emx}{\end{array}\right)}
\newcommand{\A}{\mathbb{A}}

\newcommand{\N}{\mathbb{N}}
\newcommand{\Z}{\mathbb{Z}}
\newcommand{\R}{\mathbb{R}}
\newcommand{\C}{\mathbb{C}}
\newcommand{\g}{\mathfrak{g}}
\newcommand{\cc}{{\mathcal{C}}}
\newcommand{\dd}{{\mathcal{D}}}

\newcommand{\aaa}{{\mathcal{A}}}

\newtheorem{theorem}{Theorem}[section]
\newtheorem{proposition}[theorem]{Proposition}
\newtheorem{corollary}[theorem]{Corollary}
\newtheorem{lemma}[theorem]{Lemma}
\newtheorem{defn}[theorem]{Definition}

\numberwithin{equation}{section}

\def\ni{\noindent}

\def\la{\lambda}

\def\al{\alpha}

\def\sp{\mathfrak{sp}}
\def\so{\mathfrak{so}}

\def\Ta{{\mathcal T}}

\def\H{\mathcal H}
\def\A{\mathcal A}
\def\Ca{\mathcal C}
\def\D{\mathcal D}
\def\BD{\mathcal BD}

\def\I{\mathcal I}
\def\J{\mathcal J}

\def\P{\mathcal P}

\def\A{\mathcal A}
\def\E{\mathcal E}

\def\e{\tilde e}

\def\M{\mathcal M}
\def\B{{\mathcal B}}

\def\ignore#1{\relax}

\def\Tt{\tilde T}

\def\D{{\mathcal D}}

\def\ep{\epsilon}

\def\ve{\varepsilon}

\def\Pm{{ P_{[1^m]}}}
\def\Pm1{{ P_{[1^{m+1}]}}}
\def\Pm2{{ P_{[1^{m+2}]}}}
\def\P2m{{ P_{[2,1^{m-1}]}}}
\def\Pmm1{{ P_{[1^{m-1}]}^{(m+1)}}}

\def\choose #1 #2{\begin{pmatrix}#1\\#2\end{pmatrix}}

\def\r{\tilde r}
\def\trd{\tr_{\mathcal D}}

\begin{document}

\title{On braided tensor categories of type $BCD$}

\author{Imre Tuba}
\address{Department of Mathematics, University of California, 
Santa Barbara, CA 93106-3080}
\email{ituba@math.ucsb.edu}

\author{Hans Wenzl}
\address{Department of Mathematics, University of California, 
San Diego, CA 92093-0112}
\email{hwenzl@ucsd.edu}
\thanks{H.W. was partially supported by NSF grant DMS 0070848}

\date{\today}

\keywords{tensor category, monoidal category, braiding, quantum group}

\subjclass{Primary 20F36, 20C07, 81R10;
Secondary 20H20, 16S34}

\begin{abstract}
We give a full classification of all braided semisimple
tensor categories whose Grothendieck semiring is the one of
$\rep\big(O(\infty)\big)$ (formally), $\rep\big(O(N)\big)$,
$\rep\big(Sp(N)\big)$ or of one of its associated fusion categories.
If the braiding is  not symmetric, they are completely determined by the
eigenvalues of a certain braiding morphism, and we determine precisely
which values can occur in the various cases. If the category allows
a symmetric braiding, it is essentially determined by the dimension
of the  object corresponding to the vector representation.
\end{abstract}

\maketitle

\section{Introduction}

Braided tensor categories have played a prominent role in various areas 
in recent years, such as conformal field theory, string theory, 
operator algebras and low-dimensional topology.
Important examples have been constructed in a mathematically rigorous way
using the representation theory of quantum groups, loop groups and 
Kac-Moody algebras. This naturally leads to the question of classifying
such categories. We solve this question in this paper for braided categories
associated to the representation categories of orthogonal and symplectic
groups, and various generalizations of them.

It has been shown in \cite{kazhdan-wenzl}  that any rigid semisimple
tensor category whose Grothendieck semiring is equivalent to the one of 
$\rep\big(SU(N)\big)$ must necessarily be equivalent to the
category  $\rep(U_q\sln_N)$,
with $q$ not a root of unity, up to $N$ possible choices of a twist;
here $U_q\sln_N$
is the Drinfeld-Jimbo $q$-deformation of the universal enveloping 
algebra $U\sln_N$.
The present paper 
proves a similar statement for a braided tensor category 
$\Ca$ whose Grothendieck 
semiring is isomorphic to the one
of a full orthogonal or a symplectic group. It will be convenient to
formulate the result in a slightly different way in this case:
Let $X$ be the object in $\Ca$
corresponding to the vector representation of an orthogonal or
symplectic group. It is well-known that its second tensor power
decomposes into the direct sum of three irreducible objects.
Hence the braiding morphism $c_{X,X}$ has at most three different
eigenvalues. It is easy to see that one can also define  braiding
structures for $\Ca$  by replacing $c_{X,X}$ by its inverse,
its negative or its negative inverse.
If $c_{X,X}$ has three distinct eigenvalues,
$\Ca$ is completely classified as a monoidal category
by these eigenvalues. Another set of eigenvalues
belongs to a category equivalent to $\Ca$ if and only if it can be obtained from the ones
of $c_{X,X}$ by changing the braiding structure as just mentioned
before. Moreover, we also show that the eigenvalues have to be
of the form $q, -q^{-1}$ and $r^{-1}$, or of the form 
$iq, -iq^{-1}$ and $ir^{-1}$, with $q$ not a root of unity
and with $r$ being $\pm $ a power of $q$,
where the exponent depends on the particular orthogonal or symplectic
group. Here the two possible forms of the eigenvalues correspond
to the two possible twists (in the language of \cite{kazhdan-wenzl}) 
for categories of this type. If $c_{X,X}$ has only two distinct eigenvalues, 
they are necessarily of the form $\{\pm 1\}$ or $\{\pm i\}$, and the category 
is completely determined by this and the quantity $d(X)$, which, up to a sign,
is equal to the categorical dimension of the object $X$. In particular,
we obtain two distinct families of categories whose Grothendieck semirings
are isomorphic to the one of an odd-dimensional orthogonal group,
while there is only one such family if the Grothendieck semiring
is the one of an even-dimensional orthogonal, a symplectic
or a special unitary group (see Cor \ref{catclass} for a more precise
statement).

It is easy to define a Grothendieck semiring which 
could be considered as the one of a formal group $O(\infty)$,
and one can define categories with such a Grothendieck semiring.
The methods in our paper apply similarly to classify such
categories, and we obtain essentially the same classification
as in the last paragraph. The only difference is that now
$r$ can not be $\pm$ a power of $q$, and $q$ can not be a root 
of unity. Finally, our methods also apply to fusion categories
whose Grothendieck semirings are quotients of the ones of
an orthogonal or symplectic group. Here $q$ has to be a root of 
unity and $r$ is a power of $q$, where the order of the root
of unity and the exponent depend on the given Grothendieck semiring.
We also remark that in our context the braiding condition is strong
enough that we never need to consider the full Grothendieck semiring;
it suffices to know how to tensor with the vector representation.

The method of proof in this paper is  similar to the one in 
\cite{kazhdan-wenzl}. We first give an
intrinsic description of the endomorphisms of tensor powers of an 
object $X$ corresponding
to the vector representation of an orthogonal or symplectic group 
in terms of certain representations of braid groups. From this one can
reconstruct the whole category, similarly as it was done
in \cite{kazhdan-wenzl}. In this paper, we do this following an alternate
approach due to Alain  Brugui\` eres.
Besides that, the main differences to the paper \cite{kazhdan-wenzl} are
that we have  to assume a priori that these categories 
are braided (which may not be necessary)
and that the braid representations as well as the
combinatorics involved here are more 
complicated than the ones in \cite{kazhdan-wenzl}. 

Here are the contents of this paper in more detail. We first recall
basic definitions of braided rigid tensor categories. We then present 
reconstruction techniques of \cite{kazhdan-wenzl} and from
Brugui\` eres' unpublished lecture notes \cite{bruguieres}; in particular, 
Section \ref{extending} closely follows these
notes. In Section \ref{BCD specifics},
we derive relations for the braid representations occurring in
$\End\big(X^{\otimes n}\big)$. We then study the corresponding abstract
algebras given by these relations, which depend on two parameters.
The main difficulty then is to show that these algebras map
surjectively onto $\End\big(X^{\otimes n}\big)$.  Here the crucial idea
is, as in \cite{kazhdan-wenzl}, the abstract characterization of the
trace functional on $\End\big(X^{\otimes n}\big)$ coming from the
dimension function as a so-called Markov trace. This shows that the
image has to contain at least the quotient of this algebra modulo the
annihilator ideal of the Markov trace. Rigidity is then used to prove
that the image actually has to be equal to the quotient. This result
together with the reconstruction results in Sections \ref{catreconstruction}
and \ref{extending} is then used to prove the classification result
in the last section.

\emph{Acknowledgments:} Hans Wenzl would like to thank David Kazhdan
for useful discussions. Imre Tuba would like to thank Ken Goodearl for
the same.  Both authors would like to thank Alain Brugui\` eres for
allowing them to use his unpublished notes \cite{bruguieres} in this
paper and the referee for the thorough reading and useful remarks
which improved the presentation.

\section{Definitions and notation}

We recall some basic definitions and set up notations. For more
details, we refer to \cite{maclane-cat}, \cite{freyd} for general
categorical notions, and to \cite{kassel}, \cite{turaev} for
tensor categories.
\begin{defn}
A \emph{monoidal category} $\cc$ is a category $\cc$ with a functor
$\otimes : \cc \times \cc \to \cc$ called the \emph{tensor product},
a natural isomorphism
$a$ between $\otimes \circ \left(\otimes \times \iden_{\cc}\right)$ 
and $\otimes
\circ \left(\iden_{\cc}\times \otimes\right)$ called the \emph{associativity
constraint}, satisfying the \emph{pentagon axiom},
a unit object $\1 \in \cc$ and natural isomorphisms 
$l_X : \1 \otimes X \to X$ and $r_X : X \otimes \1 \to X$ called the left 
and right \emph{unit constraints} satisfying the \emph{triangle axiom}.
\end{defn}
The pentagon axiom just states that different ways of rebracketing
the tensor product of four objects will lead to the same results,
see e.g. \cite{kassel} for a precise statement. The triangle axiom
just states that the left and right constraints are compatible with
associativity, i.e. that $(1_X\otimes l_Y)\circ a_{X,\1,Y}$ and 
$r_X\otimes 1_Y$
describe the same morphism from $(X\otimes \1)\otimes Y$ to $X\otimes Y$;
here $a_{X,\1,Y}$ is the associativity morphism $(X\otimes \1)\otimes Y
\to X\otimes (\1\otimes Y)$.
A \emph{monoidal functor} is a triple $(F, \theta, \phi)$, where $F: \cc \to \cc'$ is a functor, $\theta \in \Hom_{\cc'} \big(F(\1), \1'\big)$ is an isomorphism and $\phi$ is a natural isomorphism
\[ \phi_{X,Y} : F(X) \otimes' F(Y) \to F(X \otimes Y). \]
In order to respect the monoidal structure, $\theta$ and $\phi$ and required to satisfy certain obvious commutative diagrams. See e.g.\ \cite{kassel}, Ch.\ XI.4 for the full definition.

A monoidal category $\cc$ is called strict if $a$, $l$, and $r$ are the 
identity. That is $\left( X \otimes Y \right) \otimes Z = 
X \otimes \left( Y \otimes Z \right)$ and $\1 \otimes X = X \otimes \1 = X$ 
for any $X \in \cc$.
A theorem of Mac Lane's asserts that any monoidal category is equivalent to 
a strict one (see e.g.\ \cite{kassel}, p.\ 288). Since our interest is 
in characterizing tensor categories up to equivalence, we may and will
assume our categories to be strict monoidal for the rest of the paper.

A strict monoidal category $\cc$ is called \emph{right rigid} if every object $X \in \cc$ has a right dual object $X^* \in \cc$ and a pair of morphisms $i_X : \1 \to X \otimes X^*$ and $d_X:X^* \otimes X
\to \1$ such that the maps
\[ \xymatrix{X = \1 \otimes X \ar[rr]^-{i_X \otimes \iden_X} & & X \otimes X^* \otimes X \ar[rr]^-{\iden_X \otimes d_X} & & X \otimes \1 = X} \]
\[ \xymatrix{X^* = X^* \otimes \1 \ar[rr]^-{\iden_{X^*} \otimes i_X} & & X^* \otimes X \otimes X^* \ar[rr]^-{d_X \otimes \iden_{X^*}} & & \1 \otimes X^* = X^*} \]
are $\iden_X$ and $\iden_{X^*}$.

With this notion of duality, we also have the usual isomorphism
between $\Hom(V, W\otimes X^*)$ and $\Hom(V\otimes X, W)$ for any
objects $V, W$ in $\cc$. One checks
easily that these isomorphisms are given by
the maps $a\to (\iden_W\otimes d_X) \circ (a\otimes \iden_X)$
and
$b\to (b\otimes \iden_{X^*}) \circ (\iden_V\otimes i_X)$ 
for $a \in \Hom(V, W\otimes X^*)$ and $b \in \Hom(V\otimes X, W)$.
In particular, one obtains as a special case that $\dim \Hom(\ident,
X \otimes X^*) = \dim \End(X) = 1$ if $X$ is a simple object.
Left rigidity is defined similarly as right rigidity with the left dual ${}^*X$ of $X$ on the opposite side of $X$.

A \emph{tensor category} is an abelian category with the additional
structure of a monoidal category such that the tensor product and the
direct sum are distributive.

\begin{defn}
A \emph{$\C$-category} $\cc$ is an additive category in which the morphisms between any two objects form a finite dimensional $\C$-vector space and composition of morphisms is bilinear relative to the vector space structure. A tensor category which is also a $\C$-category will be called a \emph{$\C$-tensor category}. In this case, we will require that the categorical tensor be $\C$-bilinear.
\end{defn}

A strict monoidal category $\cc$ is called \emph{braided} if there exists a natural isomorphism $c_{X,Y}: X \otimes Y \to Y \otimes X$ called the \emph{braiding} such that:
\[ \xymatrix{
X \otimes Y \otimes Z
\ar[rr]_-{c_{X,Y \otimes Z}}
\ar[dr]_-{c_{X,Y} \otimes \iden_Z}
& & Y \otimes Z \otimes X \\
& Y \otimes X \otimes Z
\ar[ur]_-{\ \ \iden_Y \otimes c_{X,Z}} & } \]
and
\[ \xymatrix{
X \otimes Y \otimes Z
\ar[rr]_-{c_{X \otimes Y,Z}}
\ar[dr]_-{\iden_X \otimes c_{Y,Z}}
& & Z \otimes X \otimes Y \\
& X \otimes Z \otimes Y
\ar[ur]_-{c_{X,Z} \otimes \iden_Y} & } \]
are commuting diagrams. Naturality means that for any morphisms $f:X \to X'$ and $g:Y \to Y'$
\[ (g \otimes f) \circ c_{X,Y} = c_{X',Y'} \circ (f \otimes g). \]

Let $\cc$ and $\cc'$ be strict braided monoidal categories. 
A monoidal functor $(F, \theta, \phi)$ is called \emph{braided} if
it respects the braiding axioms in the sense that
\[ F(c_{X,Y})\,\phi_{X,Y} = \phi_{Y,X}\,c'_{F(X),F(Y)}. \]

A braiding is a generalization of the flip, which is the natural isomorphism $P_{A,B}:A \otimes B \to B \otimes A$ on the category of modules over the commutative ring $R$. Note that the flip is involutive, that is $P_{B,A} \circ P_{A,B} = \iden_{A \otimes B}$. This is not required for a braiding, but the property is generalized in the notion of the \emph{twist}, which is a natural isomorphism $\theta_V: V \to V$ in a braided monoidal category $\cc$ such that
\[ \theta_{X \otimes Y} = c_{Y,X} \circ c_{X,Y} \circ \left( \theta_X \otimes \theta_Y \right) \]
for all $X,\, Y \in \cc$. $\theta$ is required to be functorial in the sense that for any morphism $f:X \to Y$, $\theta_Y \circ f = f \circ \theta_X$.
A \emph{ribbon category} $\cc$ is a rigid braided monoidal category with a compatible twist, meaning:
\[ \big(\theta_X \otimes \iden_{X^*}\big) \circ i_X = \big(\iden_{X} \otimes \theta_{X^*}\big) \circ i_X. \]
In a ribbon category, right rigidity also implies left rigidity and vice versa. In fact, given the right duality morphisms $i$ and $d$,
\begin{equation}\label{rightdual}
i_X' = (\iden_{X^*} \otimes \theta_X) \circ c_{X,X^*} \circ i_X
\quad
{\rm and} \quad 
d_X' = d_X \circ c_{X,X^*}\circ (\theta_X \otimes \iden_{X^*})
\end{equation}
yield left duality morphisms which make the category left rigid. With this left duality, the left and right duals of objects and morphisms coincide.

We will also need the morphism
\begin{equation}\label{etilde}
\tilde e_X = i_X' \circ d_X = i_X \circ d_X' \in \End\big(X\otimes X^*\big).
\end{equation}
These allow us to define the \emph{categorical trace} of an endomorphism $f \in \End(X)$ as
\[ \Tr_X(f) = d_X' \circ  (f \otimes \iden_{X^*} ) \circ i_X
= d_X \circ c_{X,X^*} \circ \big(\left( \theta_X \circ f \right) \otimes \iden_{X^*}\big) \circ i_X  \in \End(\1), \]
which is easily seen to be the same as
\[ \Tr_X(f) = d_X \circ (\iden_{X^*} \otimes f) \circ i_X'
= d_X \circ \big(\iden_{X^*} \otimes \left( \theta_X \circ f \right)\big) \circ c_{X,X^*} \circ i_X  \in \End(\1), \]
using naturality of the braiding and the twist. 
Just like the usual trace of linear operators, $\Tr_Y(fg) = \Tr_X(gf)$ for any $f \in \Hom(X,Y)$ and $g \in \Hom(Y,X)$, and $\Tr_{X \otimes Y} (f \otimes g) = \Tr_X(f) \Tr_Y(g)$ for any $f \in \End(X)$ and $g \in \End(Y)$ (see 
\cite{kassel} or \cite{turaev} for a proof). 
If $f \in \End(\1)$, then $\Tr_{\smone}(f) = f$. The \emph{categorical dimension} of an object $X$ is $\dim X = \Tr_X\big(\iden_X \big)$. It is clear from the properties of the trace that $\dim X \oplus Y = \dim X + \dim Y$ and $\dim X \otimes Y = (\dim X)(\dim Y)$.

The \emph{normalized trace} $\tr_X$
on $\End(X)$  is defined by $\tr_X(f)=\Tr_X(f)/(\dim X)$. In the following
we will often just write $\Tr$, $\tr$  for the trace or normalized trace
when it is clear for which object it is defined.

We call a morphism a monomorphism or monic if its kernel is $0$ and an epimorphism or epic if its cokernel is $0$.
As is customary, we won't get hung up on abusing the language slightly and calling object $A$ a ``subobject'' of $B$ if there exists a monomorphism $A \to B$, and referring to a monomorphism in the kernel of $f$ as ``a kernel.''

\section{Categorical reconstruction}
\label{catreconstruction}

In the following we will say that a $\C$-category $\cc$ is \emph{semisimple} 
if every endomorphism ring in $\cc$ is a semisimple $\C$-algebra. 
An object $Y$ in $\cc$ is called \emph{simple} if $\End(X)=\C$.
This is a somewhat weaker definition  of semisimplicity as is usually
common, as can be seen at the following example.

\begin{defn}
A \emph{monoidal algebra} $\aaa$ is a semisimple monoidal category whose 
objects are the natural numbers with ordinary addition as the tensor product.
\end{defn}

To get an example of a monoidal algebra, let $\cc$ be a semisimple
monoidal category, and let $X$ be an object in $\cc$. 
Then the subcategory $\aaa$ whose objects are 
tensor powers of $X$ (with the obvious labeling 
$X^{\otimes n} \longleftrightarrow n\in\N$ ) is a monoidal algebra;
here we define $X^{\otimes 0}=\1$, the trivial object.
It is well-known that if one takes for $X$ the vector representation
of a classical Lie group, the only simple objects in the corresponding
monoidal algebra would be $\1$ and $X$ itself.

However, it is  well-known  that
the representation category of a classical Lie
group is essentially determined if one understands the decomposition
of tensor powers of its vector representation.
This statement will be
made precise and proved in this and the following section for 
general monoidal semisimple $\C$-tensor categories.

Let $\cc$ be a monoidal  category. In order to get  direct
sums (i.e. an additive category), we first define a 
larger category $\add \cc$ whose objects are 
finite sequences of objects from $\cc$ including the empty sequence. 
The morphisms from $\big( X_1, X_2, \ldots , X_n \big)$ to $\big( Y_1, Y_2, \ldots , Y_m \big)$ are defined by
\[ \Hom_{\add \cc} \Big(\big( X_1, X_2, \ldots , X_n \big),\,\big( Y_1, Y_2, \ldots , Y_m \big)\Big) = \bigoplus_{i,j} \Hom_{\cc} \big(X_i, Y_j\big) \]
where $\oplus$ on the right-hand side stands for the ordinary direct sum of vector spaces. If either of the two sequences is empty, the Hom space will be the $0$-vector space. We will compose morphisms, when possible, by ordinary matrix multiplication. We claim that this is an additive category with concatenation of sequences as the direct sum operation. The required direct sum system
\[ \xymatrix{
\bigl( X_1, X_2, \ldots , X_n \bigr) \ar@<.5ex>[r]^<>(.5){\iota_1} &
\bigl( X_1, X_2, \ldots , X_n, Y_1, Y_2, \ldots , Y_m \bigr) \ar@<.5ex>[l]^<>(.5){\pi_1} \ar@<-.5ex>[r]_<>(.5){\pi_2} &
\bigl( Y_1, Y_2, \ldots , Y_m \bigr) \ar@<-.5ex>[l]_<>(.5){\iota_2}} \]
is constructed the obvious way from identities in $\End\big( X_i \big)$ and $\End\big( Y_j \big)$ and zeros in the other components.

We still need to get enough subobjects.
This will be accomplished by a process called the idempotent completion 
(see \cite{freyd}, p.\ 61), which goes as follows. Starting with any category $\cc$, 
let the objects of $\idem \cc$ be the pairs $(X, p)$ where $X \in \Ob(\cc)$ and 
$p \in \End(X)$ with $p^2 = p$, that is $p$ is an idempotent. 
The morphisms in $\idem \cc$ are defined as follows
\[ \Hom_{\idem \cc} \big((X, p), (Y, q)\big) = \big\{ f \in \Hom_{\add \cc}(X, Y) \mid f = qf = fp \big\}. \]
We will say the idempotent $p$ splits if it can be factored as $p=ab$ with $a$ monic and $b$ epic. In this case, it it easy to see $ba = 1$
(identity of the source of $a$) by canceling $a$ on the left and 
$b$ on the right from $ab = p = p^2 = abab$. 
It is an easy exercise to show that idempotents split in $\idem \cc$.

Before we prove that these constructions indeed produce an abelian 
category, we will need a lemma about the existence of quasi-inverses.

\begin{lemma}
\label{quasiinverse}
Let $\cc$ be a semisimple additive $\C$-category and $f \in \Hom(X, Y)$ for some objects $X, Y$. Then there exists $g \in \Hom(Y, X)$ with $f = fgf$ and 
with $P=fg$ and $Q=gf$ idempotents in $\End(Y)$ and $\End(X)$ respectively. If $f$ is monic, then $Q=\iden_X$ and $P$ splits as $fg$. If $f$ is epic, then $P=\iden_Y$ and $Q$ splits as $gf$. 
\end{lemma}

\begin{proof} We can naturally embed $\Hom(X,Y)$ and $\Hom(Y,X)$ into
$\End(X\oplus Y)$. Hence we can consider $f$ as an element in 
$\End(X\oplus Y)$, which is semisimple. Restricting to a simple
component, it suffices to consider
 $f' \in M_n(\C)$, acting on $V = \C^n$. 
Let $V_1 = \ker(f')$ and $W_1$ such that $V_1 \oplus W_1 = V$. 
Let $W_2 = \im(f')$ and $V_2$ such that $V_2 \oplus W_2 = V$. 
Hence $f'|_{W_1} : W_1 \to W_2$ is an isomorphism. 
Let $g': W_2 \to W_1$ be the inverse of $f$. 
Extend $g'$ to $V$ by letting it act as $0$ on $V_2$. 
Doing this  for each simple component of $\End(X\oplus Y)$,
we obtain an element $\tilde g\in \End (X \oplus Y)$ such that 
$f = f\tilde gf$. Then  $g = \pi_1 \tilde g \iota_2 \in \Hom(Y, X)$
satisfies $fgf=f$, where
\[ \xymatrix{
X \ar@<.5ex>[r]^<>(.5){\iota_1} & X \oplus Y \ar@<.5ex>[l]^<>(.5){\pi_1} \ar@<-.5ex>[r]_<>(.5){\pi_2} & Y \ar@<-.5ex>[l]_<>(.5){\iota_2}} \]
is a direct sum system in $\cc$ 

That $P^2=P$ and $Q^2=Q$ is trivial. If $f$ is monic, cancel $f$ on the left from $f=fgf$ to get $\iden_Y = gf$, which makes $g$ necessarily epic, hence $P$ splits as $fg$ and similarly for $f$ epic.
\end{proof}

\begin{theorem}
Let $\cc$ be a semisimple $\C$-category. 
Then $\ab \cc = \idem \add \cc$ is a semisimple abelian category .
\end{theorem}

\begin{proof} 
The fact that $\ab \cc$ has direct sums (i.e. it is an additive
category) follows easily by applying the construction at the beginning
of this section to objects of  $\idem\add\cc$. This is left to the
reader.
To show that $\ab \cc$ is also abelian, we need to check

1.Every morphism $f \in \Hom \big((X, p), (Y, q)\big)$ must have a kernel and a cokernel. Let us construct a kernel. Let
\[ I = \big\{ g \in \End (X, p) \mid fg = 0 \big\}. \]
Clearly, $I$ is a right ideal of $\End (X, p)$, hence $I = P\, \End (X, p)$ for some idempotent $P \in \End (X, p)$ by semisimplicity. We would like to claim that $P: (X, P) \to (X, p)$ is a kernel of $f$. $P$ is monic by definition: if $P\alpha_1 = P\alpha_2$ for some $\alpha_1, \alpha_2 \in \Hom \big((Z, r), (X, P)\big)$ then
\[ \alpha_1 = P\alpha_1 = P\alpha_2 = \alpha_2. \]
That $fP = 0$ is clear. Now, suppose $fg = 0$ for some $g \in \Hom \big((Z, r), (X, p)\big)$. We will show $g$ factors through $P$. By Lemma \ref{quasiinverse}, we have $h \in \Hom \big((X,p), (Z,r)\big)$ such that $g = ghg$. Now $f(gh) = (fg)h = 0$, hence $gh \in I = P\, \End (X, p)$. Thus
\[ Pg = P(ghg) = (Pgh)g = (gh)g = g.\]
The dual construction will give a cokernel of $f$.

2. We need to show that every monomorphism is a kernel and every epimorphism is a cokernel. Let $f \in \Hom \big((X, p ), (Y, q)\big)$ be a monomorphism. Invoke Lemma \ref{quasiinverse} again to find $g \in \Hom \big((Y, q), (X, p)\big)$ such that $f = fgf$. Let $P = 1 - fg \in End (Y, q)$. Clearly, $Pf = 0$. If $h \in \Hom \big((Z, r), (Y, q)\big)$ such that $Ph = 0$, then $h = fgh$, hence $h$ factors through $f$. As $f$ is already monic, $f$ is a kernel of $P$. The dual argument shows that every epimorphism is a cokernel.

That $\ab \cc$ is semisimple is clear as its endomorphism rings are subalgebras of the endomorphism rings in $\add \cc$, which are semisimple.
\end{proof}

If $\cc$ is a monoidal category to begin with the tensor functor $\otimes$
on $\cc$ is extended to a tensor product $\otimes_{\ab \cc}$ in the 
resulting abelian category $\ab \cc$ in the obvious way as follows:
In $\add \cc$, define $\otimes_{\add \cc}$ as
\[ \big( X_1 \oplus X_2 \oplus \cdots \oplus X_n \big) \otimes_{\add \cc}  \big( Y_1 \oplus Y_2 \oplus \cdots \oplus Y_m \big)= \bigoplus_{i,j} \big( X_i \otimes_{\cc} Y_j \big) \]
on the objects and analogously on the morphisms. (Where $\oplus$ is the categorical direct sum constructed previously.) In $\ab \cc$ , define $\otimes_{\ab \cc}$ as
\[ (X, p) \otimes_{\ab \cc} (Y, q) = \big(X \otimes_{\add \cc} Y, p \otimes_{\add \cc} q\big) \]
on the objects and simply as $\otimes_{\add \cc}$ on the morphisms.

We also observe that if $\dd$ is a full subcategory of a semisimple
additive category $\cc$, then $\add \dd$ is equivalent to the additive 
subcategory generated by $\dd$ in $\cc$, that is the full subcategory 
whose objects are finite direct sums of objects of $\dd$ inside $\cc$.
We will in the following identify $\add\dd$ with that subcategory
to simplify notation.

\begin{theorem}
\label{nonmonoidal-reconstruction}
Let $\cc = (\cc, \otimes, a, \1, l, r)$ be a semisimple abelian $\C$-category and $\dd$ a full subcategory (not necessarily abelian) of $\cc$ that generates $\cc$ in the sense that every object in $\cc$ is a subobject of a direct sum of objects from $\dd$. Then there is an equivalence of abelian categories:
\[ \ab \dd \cong \cc. \]
\end{theorem}
\begin{proof}
We will construct the equivalence $F : \cc \to \ab \dd$. Let $A \in \Ob(\cc)$. For every such object, we can choose $X_1, \ldots X_n \in \Ob(\dd)$ and a monic $f: A \to X_1 \oplus \cdots \oplus X_n$ in $\cc$ by the hypothesis.
Use Lemma \ref{quasiinverse} in $\add \dd$ to find $g: X_1 \oplus \cdots \oplus X_n \to A$ such that $f = fgf$. As $fg \in \End\big( X_1 \oplus \cdots \oplus X_n \big)$ is an idempotent, we can set $F(A) = \big( X_1 \oplus \cdots \oplus X_n, fg \big)$.

Now, let $\sigma \in \Hom_{\cc}(A,B)$. As above, there exist monomorphisms $f: A \to X_1 \oplus \cdots \oplus X_n$ and $h: B \to Y_1 \oplus \cdots \oplus Y_m$ in $\cc$, and $g$ and $k$ such that $f=fgf$ and $h=hkh$. Then we already have $F(A) = \big( X_1 \oplus \cdots \oplus X_n, fg \big)$ and $F(B) = \big( Y_1 \oplus \cdots \oplus Y_m, hk \big)$. Set $F(\sigma) = 
h \sigma g$. That this is indeed in $\Hom_{\ab \dd} \Big(\big( X_1 \oplus \cdots \oplus X_n, fg \big), \big( Y_1 \oplus \cdots \oplus Y_m, hk \big)\Big)$ follows from
\[ h \sigma g = hk(h \sigma g) = (h \sigma g)fg. \]
$F$ as a map $\Hom_{\cc}(A,B) \to \Hom_{\ab \dd} \Big(\big( X_1 \oplus \cdots \oplus X_n, fg \big), \big( Y_1 \oplus \cdots \oplus Y_m, hk \big)\Big)$ in fact has an obvious inverse $G$ that takes $\phi \in \Hom_{\ab \dd} \Big(\big( X_1 \oplus \cdots \oplus X_n, fg \big), \big( Y_1 \oplus \cdots \oplus Y_m, hk \big)\Big)$ to $k \phi f$.

We have just proven that $F$ is full and faithful. It is now enough to show that each object in $\ab \dd$ is isomorphic to one in the image of $F$ (see \cite{maclane-cat}, p.\ 93) to conclude that $F$ is an equivalence.

Let $\big( Y_1 \oplus \cdots \oplus Y_m, p \big)$ be an object in $\ab \dd$. Then $p$ is an idempotent in $\End_{\cc}\big( Y_1 \oplus \cdots \oplus Y_m \big)$. In an abelian category, every morphism has a factorization into an epimorphism followed by a monomorphism (see \cite{maclane-cat}, p.\ 199). So in particular, idempotents split. Let $p$ split as $ab$ and set $A=S(a)$. Then $a: A \to Y_1 \oplus \cdots \oplus Y_m$ is a subobject, and we claim $F(A)$ is isomorphic to $\big( Y_1 \oplus \cdots \oplus Y_m, p \big)$. For suppose that in the construction of $F$ above we chose the subobject $f:A \to X_1 \oplus \cdots \oplus X_n$ and $F(A) = \big( X_1 \oplus \cdots \oplus X_n, fg \big)$. Then it is easy to verify that $ag$ is an isomorphism in
\[ \Hom_{\ab \dd} \Big(\big( X_1 \oplus \cdots \oplus X_n, fg \big), \big( Y_1 \oplus \cdots \oplus Y_m, ab \big)\Big)\]
with inverse $fb$.
\end{proof}

Note that we are making a lot of arbitrary choices in constructing this equivalence. This is to 
be expected though, as equivalences are usually not unique. Compare this with isomorphism between 
groups: one can normally find several different isomorphisms between two isomorphic groups.

In fact, a closer look at $F$ reveals that if $\cc$ is a monoidal category 
and $\dd$ is a submonoidal category, then $F$ extends to a monoidal functor. 
The proof is long and tedious, but is straightforward 
and merely an exercise in applying definitions, so we will omit it here. 
Hence $F$ is an equivalence of tensor categories and we have

\begin{theorem}
\label{reconstruction}
Let $\cc$ be a semisimple tensor category and $\dd \subseteq \cc$ a full submonoidal category. Suppose that $\dd$ generates $\cc$ in the sense that every object in $\cc$ is a subobject of a direct sum of objects from $\dd$. Then there is an equivalence of tensor categories:
\[ \ab \dd \cong \cc. \]
\end{theorem}
We will use this result in the following context:
 Let $\cc$ be a semisimple tensor category, and let $X$ be an
object in $\cc$ which generates $\cc$ in the sense that every simple
object of $\cc$ is a subobject of some tensor power of $X$. Let
$\aaa(\cc ,X)$ be the monoidal algebra generated by $X$, as described
at the beginning of this section. 
Then the  monoidal algebra $\aaa(\cc, X)$ obviously 
inherits the braiding, and it is straightforward to show that the equivalence 
in the last theorem  is an
equivalence of braided categories. Hence we obtain
\begin{corollary} \label{reconstructioncor} With the just introduced
notations we have the equivalence of braided categories
$\ab(\aaa(\cc ,X))\cong \cc$.
\end{corollary}

\section{Extending diagonals of braided monoidal algebras}\label{extending}

The results of this section have already apeared in \cite{kazhdan-wenzl}.
Here we closely follow the
presentation which was subsequently given
by Brugui\` eres  in unpublished lecture notes \cite{bruguieres} and which
has some advantages over the original one in our context.
We would like to thank Alain Brugui\` eres for allowing us to include
this material in our paper.

The precise goal of this section will be
stated after Definition \ref{defdiag}. In the following
 $\cc$ is a semisimple (not necessarily braided)
tensor category, $X$ is an object in $\cc$ and $\aaa=\aaa(\cc ,X)$
is the associated monoidal algebra, as in the last section.

\begin{defn}
A monoidal algebra $\aaa=\aaa(\cc ,X)$ is \emph{of type $N$} if
\begin{enumerate}
\item
$\Hom_{\aaa} \big( X^{\otimes m}, X^{\otimes n} \big) = 0$ unless 
$m \equiv n \mod N$.

\item
$\1$ and $X$ are simple.

\item
$\Hom_{\aaa} \big( \1, X^{\otimes N} \big) = \Hom_{\aaa} \big( X^{\otimes N}, \1 \big) = \C$.
\end{enumerate}
\end{defn}
This, for example, holds for the monoidal algebra arising from 
the vector representation in the
representation categories of $SU(N)$ and $U_q\sln_N$, and also
for orthogonal and symplectic categories with $N=2$ (see Section 
\ref{BCD specifics}).

\begin{lemma} 
\label{pi-props}
Let $\aaa$ be a monoidal algebra of type $N$.
\ 
\begin{enumerate}
\item
There exist nonzero morphisms  $\iota:\1 \to X^{\otimes N}$ and 
$\pi:X^{\otimes N} \to \1$ such that $\pi \iota = \iden_{\smone}$ and
$\iota \pi = \Pi$ is an idempotent in $\End \big(X^{\otimes N}\big)$
independent of the choices of $\iota$ and $\pi$.

\item
$\dim \Big\{ f \in \End \big(X^{\otimes N}\big) \mid f \Pi = f =\Pi f\Big\} = 1$.

\item
For any $n \in \N$, the map $\phi : \End \big(X^{\otimes n}\big) \to \End \big(X^{\otimes n+N}\big)$ which takes $f$ to $f \otimes \Pi$ is an isomorphism onto
\[\Sigma = \Big\{ g \in \End \big(X^{\otimes n+N}\big) \big|\big. \big(\iden_{X^{\otimes n}} 
\otimes \Pi\big) \, g = g \, \big(\iden_{X^{\otimes n}} \otimes \Pi\big) = g \Big\}. \]
\end{enumerate}
\end{lemma}
\begin{proof}
Let $\iota:\1 \to X^{\otimes N}$ be a nonzero morphism.
By Lemma \ref{quasiinverse} there exists a morphism $\pi:X^{\otimes N} \to \1$ such that 
$\iota \pi \iota = \iota$. It follows that $\pi \iota= \iden_{\smone}=1 \in \C$, and $\Pi = \iota \pi \in \End \big(X^{\otimes N}\big)$ is an idempotent.
This idempotent is unique as the object $\1$ appears with multiplicity $1$ in 
$X^{\otimes N}$.

The second statement is a consequence of the last statement with $n=0$.
To prove the last statement observe
that $\phi(f) \in \Sigma$ is clear from the first property of $\Pi$. 
Let $\psi: \Sigma \to \End \big(X^{\otimes n}\big)$ be defined by
\[ \psi(g) = \big(\iden_{X^{\otimes n}} \otimes \pi\big) \, g \, \big(\iden_{X^{\otimes n}} \otimes \iota\big). \]
Then it is straightforward to check that $\psi$ is the inverse of
$\phi$, which  finishes the proof of the  lemma.
\end{proof}

\begin{defn}\label{defdiag}
The \emph{diagonal} $\dd = \Delta \aaa$ of a monoidal algebra $\aaa$ is 
the monoidal algebra with
\[ \Hom_{\dd} \big( X^{\otimes m}, X^{\otimes n} \big) = 0 \text{\ if $m \neq n$} \]
and
\[ \End_{\dd} \big( X^{\otimes n} \big) = \End_{\aaa} \big( X^{\otimes n} \big). \]
\end{defn}

We will now investigate to what extent the structure of a monoidal
algebra of type $N$ can be recovered from its diagonal. 
So let $\dd$ be a braided diagonal monoidal algebra with braiding $c$,
which is the diagonal of a (not necessarily braided) monoidal
algebra $\aaa$ of type $N$. We
attach a complex number $\Theta(\aaa)$ to $\aaa$ as follows:
\begin{equation}\label{deftheta}
\Theta(\aaa) = l_X \, \big(\pi \otimes \iden_{X}\big) \, c_{1,N} \, \big(\iden_{X} \otimes \iota\big) \, r^{-1}_X \in \End(X) = \C. 
\end{equation}
In fact, since $\aaa$ is a strict category $l_X = r_X = \iden_X$. So we are free to suppress them. We will simply denote $\Theta(\aaa)$ by $\Theta$ whenever the context is clear. Observe that $\Theta$, just like $\Pi$ depends only on $\aaa$ 
and not on the particular choice of $\pi$ and $\iota$.

We will now prove some simple results for the braided diagonals of monoidal algebras $\aaa=\aaa(\cc, X)$ of type $N$.
To keep the notation from 
becoming overwhelming, we will use the simplified notation
\[ c_{m,n} = c_{X^{\otimes m}, X^{\otimes n}} \]
for the braiding.

\begin{lemma}\label{diaglemma} 
Let $\aaa$ be a monoidal algebra of type $N$. Suppose its diagonal $\dd$ has a braiding $c$. Then we have
\label{pullthrus}
\ 
\begin{enumerate}
\item
$\big(\pi \otimes \iden_X\big) \, c_{1,N} = 
\Theta \big(\iden_X \otimes \pi\big)$
and
$c_{1,N} \, \big(\iden_X \otimes \iota\big) = \Theta \big(\iota \otimes \iden_X\big)$.

\item
$\big(\iden_X \otimes \pi\big) \, c_{N,1} = \Theta \big(\pi \otimes \iden_X\big)$ and 
$c_{N,1} \, \big(\iota \otimes \iden_X\big) = \Theta \big(\iden_X \otimes \iota\big)$.

\item
$c_{N,N} \, \big(\Pi \otimes \Pi\big) = \big(\Pi \otimes \Pi\big) \, c_{N,N} = \Pi \otimes \Pi$.
\end{enumerate}
\end{lemma}

\begin{proof}
We will prove the first statement and leave the rest to the reader.
\begin{eqnarray*}
\big(\pi \otimes \iden_X\big) \, c_{1,N} & = & \big(\pi \otimes \iden_X\big) \, \big(\iota \otimes \iden_X\big) \, \big(\pi \otimes \iden_X\big) \, c_{1,N} \\
& = &\big(\pi \otimes \iden_X\big) \, c_{1,N} \, \big(\iden_X \otimes \iota\big) \, \big(\iden_X \otimes \pi\big) = \Theta \big(\iden_X \otimes \pi\big).
\end{eqnarray*}
where the first equality holds because $\pi \iota = 1$, the second because $\iota \pi = \Pi \in \End\big(X^{\otimes N}\big)$ which is in $\dd$ and $c$ is functorial on $\dd$, and the third is by the definition of $\Theta$. The second part of the first statement goes similarly.
\end{proof}

Let $\aaa$ and $\aaa'$ be two monoidal algebras of type $N$ with
braided diagonals. We say that $\aaa$ and $\aaa'$ are 
extensions
of the diagonal $\dd=\Delta\aaa$ if there is an equivalence $\Psi$
between $\dd$ and the diagonal $\dd'$ of $\aaa'$ as braided
categories
such that $\Psi(X^{\otimes n})=(X')^{\otimes n}$ for all $n \in \N$.
We say that the extensions $\aaa$ and $\aaa'$ of $\dd$ are
\emph{diagonally
equivalent} if $\Psi$ can be extended to
an equivalence $\Phi : \aaa \to \aaa'$ of
monoidal algebras.  

We are going to show that $\Theta(\aaa)$ is an invariant
under diagonal equivalence. 

\begin{proposition}
\label{Theta-well-defined}
\ 
\begin{enumerate}
\item
 Let $\aaa$ and $\aaa'$ be monoidal algebras of type $N$ and $\Phi : \aaa \to \aaa'$ a diagonal equivalence. Then $\Theta(\aaa) = \Theta(\aaa')$.

\item
$(\Theta(\aaa))^N=1$.
\end{enumerate}
\end{proposition}

\begin{proof} 
Since $\Phi$ is a monoidal functor $\aaa \to \aaa'$, 
it comes equipped with the isomorphism $\theta: \Phi(\1) \to \1'$
and the natural isomorphism
\[ \phi_{i,j} : \Phi(X^{\otimes i}) \otimes' \Phi(X^{\otimes j}) 
\to \Phi(X^{\otimes i+j}) \]
compatible with the action of $\Phi$ on morphisms (see e.g.\ \cite{kassel}, Ch.\ XI.4). This means,
in particular, that we have for any morphisms $f: X^{\otimes i}\to
X^{\otimes r}$ and $g: X^{\otimes j} \to X^{\otimes s}$ 
\[ \Phi(f \otimes g)=\phi_{r,s}^{-1} \circ \big(\Phi(f) \otimes' \Phi(g)\big)
\circ \phi_{i,j} \]
and compatibility with the braiding means that
\[ c'_{i,j}=\phi_{j,i}^{-1} \circ \Phi(c_{i,j}) \circ \phi_{i,j}.\]
Moreover, compatibility with the left and right unit constraints
translates into the identities
\[ \Phi(l_X) \circ \phi_{0,1} = l'_{X'} \circ (\theta \otimes \iden_{X'}) \qquad \text{and} \qquad \Phi(r_X) \circ \phi_{1,0} = r'_{X'} \circ (\iden_{X'} \otimes \theta). \]
But monoidal algebras are strict monoidal categories, so the unit constraints are identities. Using the bilinearity of the tensor product and the naturality of the unit constraints we obtain
\[ \phi_{0,1} = \theta \otimes \iden_{X'} = \iden_{\smone} \otimes \, \theta\iden_{X'} = \theta\iden_{X'} \qquad \text{and} \qquad \phi_{1,0} = \iden_{X'} \otimes \, \theta = \theta\iden_{X'} \otimes \iden_{\smone} = \theta\iden_{X'}, \]
and thus $\phi_{0,1} = \phi_{1,0}$.
Now observe
\[ \Phi(\pi_{\aaa}) \, \Phi(\iota_{\aaa}) = \Phi(\pi_{\aaa} \iota_{\aaa}) = \Phi(\iden_{\smone}) = \iden_{\smone}. \]
Hence we can and will choose $\pi_{\aaa'} = \Phi\left( \pi_{\aaa} \right)$ 
and $\iota_{\aaa'} = \Phi\left( \iota_{\aaa} \right)$. 
As we pointed out, $\Theta(\aaa')$ is independent of the particular choice 
of $\pi_{\aaa'}$ and  $\iota_{\aaa'}$. Using this and the identities
above, we obtain
\begin{eqnarray*}
\Theta(\aaa') & = & \big(\pi_{\aaa'} \otimes' \iden_{X'}\big) \, c'_{1,N} \, \big(\iden_{X'} \otimes' \iota_{\aaa'}\big) \\
& = & \big(\phi^{-1}_{0,1} \circ \Phi(\pi_{\aaa} \otimes \iden_{X}) \circ \phi_{N,1}\big) 
\big(\phi_{N,1}^{-1}\,\Phi(c_{1,N})\,\phi_{1,N} \big)
\big(\phi^{-1}_{1,N}  \circ \Phi(\iden_{X} \otimes \iota_{\aaa}) \circ \phi_{1,0} \big) \\
& = & \phi^{-1}_{0,1} \, \Phi\big((\pi_{\aaa} \otimes \iden_{X})\,c_{1,N}\,(\iden_{X} \otimes \iota_{\aaa})\big) \, \phi_{1,0} = \phi^{-1}_{0,1} \, \Phi\big(\Theta(\aaa)\big) \, \phi_{1,0} \\
& = & \phi^{-1}_{0,1} \, \Theta(\aaa) \, \phi_{1,0} = \Theta(\aaa),
\end{eqnarray*}
where $\Phi\big(\Theta(\aaa)\big) = \Theta(\aaa)$ because $\End_{\aaa}(X) = \End_{\aaa'}(X) = \C$ and $\Phi(\iden_X) = \iden_{X'}$.

To prove the second statement, observe that 
$c_{n,N}(1_{X^{\otimes n}}\otimes \iota)
=\Theta^n (\iota\otimes 1_{X^{\otimes n}})$;
this follows from Lemma \ref{pullthrus}(a) by induction on $n$,
using $c_{n,N}=(c_{1,N}\otimes 1_{X^{\otimes n-1}})(1_X\otimes c_{n-1,N})$.
Hence we obtain, using  Lemma \ref{pullthrus}(c),
\begin{eqnarray*}
\Pi\otimes \Pi & = & c_{N,N}(1_{X^{\otimes n}}\otimes \iota)
\iota (\pi\otimes\pi) \\
& = & \Theta^N (\iota \otimes 1_{X^{\otimes n}})
\iota (\pi\otimes\pi) \\
& = &  \Theta^N (\iota \otimes \iota)(\pi\otimes\pi) = 
\Theta^N (\Pi\otimes \Pi).
\end{eqnarray*}
\end{proof}
\ignore{\begin{eqnarray*}
\lefteqn{\Theta \iden_X \otimes \big((\pi \otimes \iden_{X^{\otimes n}})\, c_{n,N}\, (\iden_{X^{\otimes n}} \otimes \iota)\big) =} \\
& = & \Theta \, (\iden_X \otimes \pi \otimes \iden_{X^{\otimes n}})\, (\iden_X \otimes c_{n,N})\, (\iden_{X^{\otimes n+1}} \otimes \iota) \\
& = & (\pi \otimes \iden_{X^{\otimes n+1}})\, (c_{1,N} \otimes \iden_{X^{\otimes n}})(\iden_X \otimes c_{n,N})\, (\iden_{X^{\otimes n+1}} \otimes \iota) \\
& = & (\pi \otimes \iden_{X^{\otimes n+1}})\, c_{n+1,N}\, (\iden_{X^{\otimes n+1}} \otimes \iota),
\end{eqnarray*}
where the second equality is by Lemma \ref{pullthrus}(a). Now use this and induction as follows
\begin{eqnarray*}
\Theta^N & = & \Theta^N \, \pi \iota = \pi \, (\Theta \iden_X)^{\otimes N} \, \iota = \pi\, (\pi \otimes \iden_{X^{\otimes N}})\, c_{N,N}\, (\iden_{X^{\otimes N}} \otimes \iota)\, \iota \\
& = & (\pi \otimes \pi) \, c_{N,N} \, (\iota \otimes \iota) = (\pi \otimes \pi) \, c_{N,N} \, (\iota \otimes \iota) \, (\pi \iota \otimes \pi \iota) \\
& = & (\pi \otimes \pi) \, c_{N,N} \, (\Pi \otimes \Pi) \, (\iota \otimes \iota) = (\pi \otimes \pi) \, (\Pi \otimes \Pi) \, (\iota \otimes \iota) = 1.
\end{eqnarray*}}

\begin{proposition}
\label{Theta-injective}
Let $\aaa$ and $\aaa'$ be monoidal algebras of type $N$ which are extensions
of a given diagonal algebra $\dd$.
If  $\Theta(\aaa) = \Theta(\aaa')$, then $\aaa$ and $\aaa'$ 
are diagonally equivalent. 
\end{proposition}

\begin{proof}
Choose $\iota_{\aaa}$, $\iota_{\aaa'}$, $\pi_{\aaa}$, and $\pi_{\aaa'}$
which satisfy the conditions of the morphisms $\iota$ and $\pi$
in Lemma \ref{pi-props} for $\aaa$ and $\aaa'$. We will construct an 
equivalence $\Phi: \aaa \to \aaa'$ of monoidal algebras extending
the equivalence $\Psi$ between their diagonals.
Define $\Phi\mid_{\dd} = \iden_{\dd}$,
$\Phi\left( \iota_{\aaa} \right)= \iota_{\aaa'}$ and 
$\Phi\left( \pi_{\aaa} \right) = \pi_{\aaa'}$.
This will ensure uniqueness of a functor $\Phi$.

If $m \equiv n \mod N$, let $p, \alpha, \beta \in \N$ such that $p = m + \alpha N = n + \beta N$. The idea is to pad $f$ with $\iota$'s on the right and $\pi$'s on the left so that the result is in $\End\big(X^{\otimes p}\big)$. Let
\begin{equation}\label{fp}
f_p = \big(\iden_{X^{\otimes n}} \otimes \iota^{\otimes \beta}\big) \, f \, 
\big( \iden_{X^{\otimes m}} \otimes \pi^{\otimes \alpha}\big) \in 
\End\big(X^{\otimes p}\big). 
\end{equation}
Note that $f_p$ is a morphism in $\Delta \aaa$. 
Multiplying the last equation by 
$\big(\iden_{X^{\otimes n}} \otimes \pi^{\otimes \beta}\big)$
from the left and by 
$\big( \iden_{X^{\otimes m}} \otimes \iota^{\otimes \alpha}\big) $ from
the right, we obtain
\begin{equation}\label{pf}
f = \big(\iden_{X^{\otimes n}} \otimes \pi^{\otimes \beta}\big) \, f_p \, 
\big( \iden_{X^{\otimes m}} \otimes \iota^{\otimes \alpha}\big).
\end{equation}
As $f_p$ is a morphism in  $\Delta(\aaa)$, we can define
\begin{equation}
\Phi(f)=\big(\iden_{X^{\otimes n}} \otimes \pi_{\aaa'}^{\otimes \beta}\big)
\, \Psi(f_p)\, 
\big( \iden_{X^{\otimes m}} \otimes \iota_{\aaa'}^{\otimes \alpha}\big).
\end{equation}
It is easy to check that
$\Phi(f)$  does not depend on the choice of $p$.
We still need to make sure that $\Phi$ is well-behaved with respect 
to the tensor product. Let $f \in \Hom_{\aaa} \big( X^{\otimes m}, X^{\otimes n} \big)$ and $\alpha$, $\beta$, $p$ such that $p = m + \alpha N = n + \beta N$. Let $g \in \Hom_{\aaa} \big( X^{\otimes m'}, X^{\otimes n'} \big)$ and $\alpha'$, $\beta'$, $p'$ such that $p' = m + \alpha' N = n + \beta' N$.
Then
\ignore{\begin{eqnarray*}
\Phi(f \otimes g) & = & \Big(\iden_{X^{\otimes n+n'}} \otimes \pi_{\aaa'}^{\otimes \beta+\beta'}\Big) \; (f \otimes g)_{p+p'} \; \Big(\iden_{X^{\otimes m+m'}} \otimes \iota_{\aaa'}^{\otimes \alpha+\alpha'}\Big).
\end{eqnarray*}
(Remember this is independent of the choice of $p$ and $p'$ anyway.)}
\begin{eqnarray*}
\lefteqn{\Phi(f) \otimes \Phi(g) =} \\
& = & \Big(\big(\iden_{X^{\otimes n}} \otimes \pi_{\aaa'}^{\otimes \beta}\big) \, \Psi(f_p) \, \big(\iden_{X^{\otimes m}} \otimes \iota_{\aaa'}^{\otimes \alpha}\big)\Big) \otimes \Big(\big(\iden_{X^{\otimes n'}} \otimes \pi_{\aaa'}^{\otimes \beta'}\big) \, \Psi(g_{p'}) \, \big(\iden_{X^{\otimes m'}} \otimes \iota_{\aaa'}^{\otimes \alpha'}\big)\Big) \\
& = & \big(\iden_{X^{\otimes n}} \otimes \pi_{\aaa'}^{\otimes \beta} \otimes \iden_{X^{\otimes n'}} \otimes \pi_{\aaa'}^{\otimes \beta'}\big) \; \Psi\big(f_p \otimes g_{p'}\big) \; \big(\iden_{X^{\otimes m}} \otimes \iota_{\aaa'}^{\otimes \alpha} \otimes \iden_{X^{\otimes m'}} \otimes \iota_{\aaa'}^{\otimes \alpha'}\big).
\end{eqnarray*}
Now use Lemma \ref{pullthrus} to move all the $\iota$'s and $\pi$'s to the right in this last expression (remember to do so in $f_p \otimes g_{p'}$), and observe that all the $\Theta$'s and $\Theta^{-1}$'s magically cancel. 
It is now clear that the expression we obtain is equal to $\Phi(f\otimes g)$.
We can construct $\Phi^{-1} : \aaa' \to \aaa$ in the analogous way, which shows that $\Phi$ is indeed an equivalence of monoidal algebras.
\end{proof}

It follows from the last two propositions that there
are at most $N$ monoidal algebras of type $N$ with the same diagonal.
Before proving their  existence, we need to determine the compatibility of
their braidings. 

\begin{proposition}
\label{braidedifftheta=1} Let $c$ be a braiding on $\dd$. Then
$c$ extends to a braiding on $\aaa$ if and only if $\Theta = 1$.
\end{proposition}
\begin{proof}
$\implies$: This is clear by functoriality.

$\Longleftarrow$: As $c$ is a braiding on $\dd$, it already satisfies most of the braiding axioms on $\aaa$ as well, except possibly functoriality. So all we have to prove is functoriality.

Now, let $f \in \Hom_{\aaa}\big(X^{\otimes m}, X^{\otimes n}\big)$. We will show $c_{1,n} \, \big( \iden_X \otimes f \big) = \big( f \otimes \iden_X \big) \, c_{1,m}$. If $m \not\equiv n \mod N$, then $f=0$ and the statement is obvious. 
Let $f \in \Hom_{\aaa} \big( X^{\otimes m}, X^{\otimes n} \big)$ and $\alpha$, $\beta$, $p$ as usual $p = m + \alpha N = n + \beta N$. Let $f_p$ be as in Eq \ref{fp}.

It follows from Lemma \ref{pullthrus} (with $\Theta=1$),
the definition of $c_{n,m}$ and from $ n+\beta N=p=  m + \alpha N $
that
\[c_{1,n} \big(\iden_{X^{\otimes n+1}} \otimes \pi^{\otimes \beta}\big) =
 \big(\iden_{X^{\otimes n}} \otimes \pi^{\otimes \beta} 
\otimes \iden_X\big) \, c_{1, p}.\]
\[c_{1, p} \big(\iden_{X^{\otimes m+1}} \otimes 
\iota^{\otimes \alpha}\big) =
 \big(\iden_{X^{\otimes m}} \otimes \iota^{\otimes \alpha} \otimes \iden_X\big) \, c_{1,m}.\]
Using this and Eq \ref{pf} we obtain
\begin{eqnarray*}
\lefteqn{c_{1,n} \big(\iden_X \otimes f\big) =} \\
& = & c_{1,n} \big(\iden_{X^{\otimes n+1}} \otimes \pi^{\otimes \beta}\big) \, \big(\iden_X \otimes f_p\big) \, \big(\iden_{X^{\otimes m+1}} \otimes \iota^{\otimes \alpha}\big)\\
& = & \big(\iden_{X^{\otimes n}} \otimes \pi^{\otimes \beta} \otimes \iden_X\big) \, c_{1, p} \, \big(\iden_X \otimes f_p\big) \big(\iden_{X^{\otimes m+1}} \otimes \iota^{\otimes \alpha}\big) \\
& = & \big(\iden_{X^{\otimes n}} \otimes \pi^{\otimes \beta} \otimes \iden_X\big) \big(f_p \otimes \iden_X\big) \, c_{1, p} \, \big(\iden_{X^{\otimes m+1}} \otimes \iota^{\otimes \alpha}\big)\\
& = & \big(f \otimes \iden_X\big) \, c_{1,m}.
\end{eqnarray*}
For $g \in \Hom_{\aaa}\big(X^{\otimes m'}, X^{\otimes n'}\big)$, a similar computation proves $c_{n',1} \, \big(g \otimes \iden_X\big) = \big(\iden_X \otimes g\big) \, c_{m',1}$. Now we use induction to conclude
\begin{eqnarray*}
\lefteqn{c_{n',n} (g \otimes f) =} \\
& = & c_{n',n} \, \big(\iden_{X^{\otimes n'}} \otimes f\big) \, \big(g \otimes \iden_{X^{\otimes m}}\big) = \\
& = & \big(f \otimes \iden_{X^{\otimes n'}}\big) \, c_{n',m} \, \big(g \otimes \iden_{X^{\otimes m}}\big) \\
& = & \big(f \otimes \iden_{X^{\otimes n'}}\big) \, \big(\iden_{X^{\otimes m}} \otimes g\big) \, c_{m',m} = (f \otimes g) \, c_{m',m}.
\end{eqnarray*}
\end{proof}

We can now prove the main result of this section. It first appeared in
\cite{kazhdan-wenzl}, with the presentation
in this section following the notes \cite{bruguieres} by Brugui\` eres.

\begin{theorem}
\label{Theta-surjective}
Let $\dd$ be the diagonal of a braided monoidal algebra of type $N$.
Then there exist exactly $N$ monoidal algebras $\aaa$ such that
$\dd=\Delta(\aaa)$ up to diagonal equivalence, one for each possible
value of $\Theta(\aaa)$.
\end{theorem}

\begin{proof}
In view of our previous results, it suffices
to construct a monoidal algebra $\aaa$ of type $N$ such that
$\Theta(\aaa)=\mu$  for each given $N$-th root of unity $\mu$.
Choose $\tau$ such that $\tau^N = 1/\mu$. Let $c'_{m,n} = \tau^{mn} c_{m,n}$. It is easy to see that this is still a braiding on $\dd$.
Denote the objects of $\aaa$ by $X^{\otimes n}$ as before. Let
\[ \Hom_{\aaa} \big( X^{\otimes m}, X^{\otimes n} \big) = 0 \quad \text{if}\  m \not \equiv n \mod N,\]
otherwise let $p = m + \alpha N = n + \beta N$ ($\alpha, \beta \in \N$)  and define
\[ A_m^n (p) = \Big\{ f \in \End_{\dd}\big( X^{\otimes p} \big) \mid \big( \iden_{X^{\otimes n}} \otimes \Pi^{\otimes \beta} \big) f = f = f \big(\iden_{X^{\otimes m}} \otimes \Pi^{\otimes \alpha}\big) \Big\}\]
 Let $A_m^n = A_m^n (p)$ with $p = \max(m,n)$. By the 3rd property of $\Pi$, we know that the map $\phi:f \mapsto f \otimes \Pi$ is an injection $\End_{\dd}\left( X^{\otimes p} \right) \to \End_{\dd}\left( X^{\otimes p+N} \right)$. Observe that the restriction of $\phi$ to $A_m^n (p)$ has exactly $A_m^n (p+N)$ for its image in $\End_{\dd}\left( X^{\otimes p+N} \right)$. Hence tensoring repeatedly on the right by $\Pi$ gives us a chain of isomorphisms 
\[ A_m^n = A_m^n (p) \cong A_m^n (p+N) \cong A_m^n (p+2N) \cong \ldots .\]
Let $\phi_P: A^n_m\to A_m^n (P)$ be the induced isomorphism, with $P\equiv p$
mod $N$. Set
\[ \Hom_{\aaa} \left( X^{\otimes m}, X^{\otimes n} \right) = A_m^n\cong
A_m^n(p)\cong A_m^n (p+N) \cong \ldots .  \]
In the following we will use  these isomorphisms to define composition
and tensor products for morphisms in $\aaa$.
 Let $g \in \Hom_{\aaa} \left( X^{\otimes k}, X^{\otimes m} \right)$ and $f \in \Hom_{\aaa} \left( X^{\otimes m}, X^{\otimes n} \right)$ with $k \equiv m \equiv n \mod N$. Choose some $P \geq \max(k,m,n)$ with
\[ P = k + \alpha N = m + \beta N = n + \gamma N \qquad \alpha, \beta, \gamma \in \N. \]
Then we define $f\circ g$ by
\[ f\circ g = \phi_P^{-1}(\phi_P(f)\circ \phi_P(g)) \]
where the three $\phi_P$'s are three different maps and are to be understood in the appropriate context.

It is easy to see that this definition is independent of the
choice of $P$. As the actual composing of maps always happens inside some 
$\End_{\dd}\big(X^{\otimes P}\big)$, associativity of this composition 
law is inherited from $\dd$.

We need to define a tensor product on this category. 
Let $f \in \Hom_{\aaa} \big(X^{\otimes m}, X^{\otimes n}\big)$ and 
$g \in \Hom_{\aaa} \big(X^{\otimes m'}, X^{\otimes n'}\big)$. 
Find $p$ and $p'$ such that $f \in A_m^n (p)$ and $g \in A_{m'}^{n'} (p')$ and
\[ p = m + \alpha N = n + \beta N \text{\ and\ } p' = m' + \alpha' N = n' + \beta' N. \]
Then
\[ \big(\iden_{X^{\otimes n}} \otimes {c'}^{-1}_{n', \beta N} \otimes \iden_{X^{\otimes \beta' N}}\big) \; \big( f \otimes_{\dd} g \big) \; 
\big(\iden_{X^{\otimes m}} \otimes c'_{m', \alpha N} \otimes \iden_{X^{\otimes \alpha' N}}\big) \]
is in $A_{m+m'}^{n+n'} (p+p')$. Applying $\phi_{p+p'}^{-1}$
to it gives us the desired morphism $f \otimes_{\aaa} g\in A_{m+m'}^{n+n'} $.
That this is strictly associative follows from the strictness of the tensor product in $\dd$ and the braiding axioms.

For $\aaa$ to be a monoidal algebra, it also needs to be a semisimple category, but that is obvious as the endomorphism rings of $\aaa$ all come from $\dd$, which is already a monoidal algebra, hence semisimple.
As
\[ \Hom_{\aaa} \big(\1, X^{\otimes N}\big) \cong A_0^N(N) = \Big\{ f \in \End_{\dd}\big(X^{\otimes N}\big) \mid \Pi f = f = f \Pi \Big\} = \C\]
and similarly for $\Hom_{\aaa} \big(X^{\otimes N}, \1\big)$, $\aaa$ satisfies all of the conditions for being a monoidal algebra of type $N$.
For $\iota$ and $\pi$ in $\aaa$, choose $\Pi$ considered as an
element in $A_0^N (N)$ and as an element in $ A_N^0 (N)$ respectively.
Then $\Pi(\aaa) = \pi \iota = \Pi^2 = \Pi$ by the 1st property of $\Pi$.
We can now verify
\begin{eqnarray*}
\Theta(\aaa) & = & \big(\pi \otimes_{\aaa} \iden_{X}\big) \; c_{1,N} \; \big(\iden_{X} \otimes_{\aaa} \iota\big) = \big(\underbrace{\Pi}_{\in A_N^0} \otimes_{\aaa} \iden_{X}\big) \; c_{1,N} \; \big(\iden_{X} \otimes_{\aaa} \underbrace{\Pi}_{\in A_0^N}\big) \\
& = & {c'}^{-1}_{1,N} \; \big(\Pi \otimes_{\dd} \iden_X\big) \; c_{1,N} \; \big(\iden_{X} \otimes_{\dd} \Pi\big) \\
& = & \tau^{-N} \; c^{-1}_{1,N} \; c_{1,N} \; \big(\underbrace{\iden_{X} \otimes_{\dd} \Pi}_{\in A_{N+1}^1(N+1)}\big) \circ \big(\underbrace{\iden_{X} \otimes_{\dd} \Pi}_{\in A_1^{N+1}(N+1)}\big) = \mu\, \underbrace{\iden_{X}}_{A_1^1} \in \End_{\aaa} (X).
\end{eqnarray*}
As $\End_{\aaa} (X) = \End_{\dd} (X) = \C$, this is just the number $\mu$.
We have just proven the existence of a monoidal algebra $\aaa$ with diagonal $\dd$ and $\Pi(\aaa) = \Pi$, and with $\Theta(\aaa)=\mu$.
\end{proof}

As we observed in Proposition \ref{braidedifftheta=1}, the braiding
$c$ on $\dd$ extends to a braiding on $\aaa$ if and only if 
$\Theta(\aaa) = 1$. If $\Theta\neq 1$ we
use the braiding $c'$ instead of $c$ as in the previous proof, which does
change $\Theta$ to 1. So the braiding $c'$ can be extended
to a braiding on $\aaa$ also in that case. It follows that all possible $N$ 
extensions $\aaa$ of $\dd$ can be given the structure of a braided
category. We have shown

\begin{corollary}\label{monclass} 
A fixed braiding of $\dd$ extends to a braiding  of only
one of the $N$ possible monoidal algebras of which it can be the diagonal.
However, for a given other monoidal algebra
$\aaa$ we can always find a braiding of $\dd$ which does extend to a braiding
of $\aaa$.
\end{corollary}

\section{Rigid Categories}\label{rigidcat}

We collect and (re)prove a number of basic results
about rigidity in braided categories which are probably well-known
to experts. This will be done in the context of ribbon tensor categories,
so we need not worry about left- or right-rigidity.

\begin{lemma}\label{rigidity}
Let $\Ca$ be a rigid semisimple ribbon tensor category. Then
any simple object has nonzero dimension. In particular,
the bilinear form $\langle a,b\rangle = \tr(a\circ b)$ on
$\End (Z)$ is nondegenerate for any object $Z$ in $\Ca$.
\end{lemma}

\begin{proof}  Let $X$ be a simple object in $\Ca$, with dual object $Y$.
Let $i_X:\ident\to X\otimes Y$, $d_X:Y\otimes X\to \ident$,
$i'_X:\ident\to Y\otimes X$, and $d'_X:X\otimes Y\to \ident$ be the
corresponding left and right rigidity morphisms. 
As $X$ is simple, the object $\ident$ appears with multiplicity
one in $X\otimes Y$. Let $\Pi$ be the unique projection onto it.
If $\dim X=0$, then $(i_X\circ d_X')^2=0$.
Hence the morphism $i_X\circ d_X'$ is a 
nilpotent multiple of $\Pi$, and therefore it must be equal to 0.
But this would contradict the rigidity axiom as follows:
$$0=\big[\iden_X\otimes d_X\big] \circ \big[(i_X\circ d_X')\otimes \iden_X\big]
\circ \big[\iden_X\otimes  i_X'\big] =$$
$$= \big[(\iden_X\otimes d_X) \circ (i_X\otimes \iden_X)\big] \circ \big[(d_X'\otimes \iden_X) \circ (\iden_X\otimes  i_X')\big] = \iden_X \circ \iden_X = \iden_X,$$
a contradiction (here the second equality follows from the rigidity axiom
and from \cite{kassel}, Prop. XIV.3.5).
\end{proof}

It will also be convenient to define partial trace operations,
which are also known under the names contractions or conditional expectations.
Let $X$ and $V$ be objects in $\Ca$.
We define the map $\ve_V$ from $\End(V\otimes X)$ onto $\End(V)$ by
\begin{equation}\label{conddef}
\ve_V(b)=\frac{1}{\dim X}
(\iden_V\otimes d'_X)\circ(b\otimes \iden_Y)\circ
(\iden_V\otimes i_X).
\end{equation}
We have the following results:

\begin{lemma}\label{condlemma} Let $b\in \End(V\otimes X)$ and let
$p=1/(\dim X)\tilde e$ be the projection onto
$\1\subset X\otimes Y$. Then
\begin{enumerate}
\item
$\tr_{V\otimes X}(b)=\tr_V\big(\ve_V(b)\big)$; in particular, if $V$ is simple
then $\ve_V(b)=\tr_{V \otimes X}(b)1$.

\item
$(1_V\otimes p)\circ(b\otimes 1_Y)
\circ(1_V\otimes p)=\ve_V(b)\otimes p$.
\end{enumerate}
\end{lemma}

\begin{proof} These statements are easy consequences from the definitions
(see also e.g. \cite{orellana-wenzl-spin-groups}, Prop. 1.4).
\end{proof}

We shall need the results of the last lemma in the following setting.
Let $m\in \End\big(X^{\otimes 2}\big)$. Then we define the morphism
$m_i\in {\rm End\ }(X^{\otimes k})$ by
$$m_i\ =\ 1_{i-1}\otimes m\otimes 1_{k-1-i},$$
where $1_r$ is the identity morphism of $X^{\otimes r}$. Then we have
the following (see also e.g. \cite{orellana-wenzl-spin-groups}, Prop. 1.4)

\begin{corollary}\label{condcor}
\ 
\begin{enumerate}
\item (Markov property) If $a\in\End\big(X^{\otimes n}\big)$, 
then $\tr\big((a\otimes 1)\circ m_n\big)=\tr(a)\tr(m)$.

\item Assume that $X$ is a self-dual object (see below) and that
$X^{\otimes 2}\cong \oplus_{j=1}^d X_{\mu_j}$, and assume that
we can write $1_X=\sum_jp_{\mu_j}$ as a sum of commuting projections 
$p_{\mu_j}\in \End\big(X^{\otimes 2}\big)$ such that $\im(p_{\mu_j})\cong X_{\mu_j}$. Then
$p_2(p_{\mu_j}\otimes 1)p_2=\frac{\dim X_{\mu_j}}{(\dim X)^2}p_2$.
\end{enumerate}
\end{corollary}

\subsection{Self-dual objects}\label{seldual}
Let $\Ca$ be a semisimple ribbon
tensor category, and let $X$ be an object in $\Ca$ which is isomorphic
to its dual. Hence  we have $i=i_X: \ident \to X^{\otimes 2}$,
$d=d_X: X^{\otimes 2}\to \ident$, $i'=i'_X: \ident \to X^{\otimes 2}$, and
$d'=d'_X: X^{\otimes 2}\to \ident$ satisfying the left and right rigidity axioms.
In the following we will denote
the braiding morphism $c_{X,X}\in \End\big(X^{\otimes 2}\big)$ just by $c$,
and $i\circ d'$ by $\tilde e$.
The morphisms $i_1$ and $i_2$ are defined by
$$i_1=i\otimes 1_1: X\cong X\otimes \ident\to X^{\otimes 3}
\quad {\rm and}\quad i_2= 1_1\otimes i : \ident\otimes X\to X^{\otimes 3},$$
with $d_1$ and $d_2$ being morphisms from $X^{\otimes 3}$ to $X$ defined
similarly.

\begin{lemma}\label{tracepreli} Let $X$ be a simple self-dual object with
dimension $\dim X$ and let
$\r\in F$ be the scalar such that $\theta_X=\r 1_X$. Then there exists
$\al\in\{\pm 1\}$ such that 
$c\circ i=\al \r^{-1}i$,
$\tr(c)=\r/(\dim X)$ and $\tr(\e)=1/(\dim X)$ for the normalized categorical
trace $\tr$ for $\End \big(X^{\otimes 2}\big)$.
\end{lemma}

\begin{proof} By definition, $\dim X=d'\circ i=\Tr(d'\circ i)=\Tr(i\circ d')=\Tr(\e)$,
which implies the statement for $\tr(\e)$.
As $\theta_{\ident} = 1_{\ident}$, it follows that 
$$i=\theta_{X\otimes X}\circ i= c\circ c \circ (\theta_X\otimes \theta_X)
\circ i = \r^2\ c\circ c\circ i.$$
As $ c\circ i$ is a multiple of $i$, the first claim follows.
This also implies that $i'=\al i$ and $d'=\al d$.
Using the braiding axioms, we obtain the identity
$c_1\circ c_2\circ i_1=i_2$; after multiplying by $d_1'$ from the right,
we obtain the equality $c_1\circ c_2\circ (i_1\circ d_1')=\al (i_2\circ d_2')
\circ  (i_1\circ d_1') $ 
in $\End\big(X^{\otimes 3}\big)$.
Using the trace property and the Markov property, we obtain
$\tr\big(c_1 \circ c_2 \circ(i_1 \circ d_1')\big)=\alpha \r^{-1}\tr(c)/(\dim X)$,
which has to be equal to $\alpha\big(\tr(i \circ d')\big)^2=\alpha/(\dim X)^2$.
The claim follows from this.
\end{proof}

The following lemma corrects the statement of Lemma 3.2 in \cite{Tuba-Wenzl-Braid};
the proof there would have been sufficient for the purposes in that
paper  and also  for this paper.

\begin{lemma}\label{irred}
The algebra generated by $\End\big(X^{\otimes 2}\big) \otimes 1$ and
by $\e_2$ acts irreducibly on the space
$\Hom\big(X,X^{\otimes 3}\big)$, via
concatenation.
\end{lemma}

\begin{proof} We use the notations as in Corollary \ref{condcor},(b),
with $p_{\mu,1}=p_\mu\otimes 1$.
Observe that 
$$(p_{\mu,1}\circ\e_2\circ p_{\nu,1})\circ(p_{\kappa,1}\circ\e_2
\circ p_{\gamma,1})=
\delta_{\nu,\kappa} (\dim X) \tr(p_\nu) (p_{\mu,1}\circ\e_2\circ p_{\gamma,1}).$$
Hence  the set $\{ p_{\mu_i,1}\circ\e_2\circ p_{\mu_j,1},\ i,j=1,2,\ ...\ d\}$
spans a full $d \times d$ matrix algebra. It
obviously does not act trivially on $\Hom\big(X,X^{\otimes 3}\big)$.
As $\dim \Hom\big(X,X^{\otimes 3}\big)=\dim \End\big(X^{\otimes 2}\big)=d$, 
the claim follows.
\end{proof}

\section{Categories of orthogonal or symplectic type}
\label{BCD specifics}

\subsection{Combinatorial data}\label{combdat}
We fix some notations for the representation category
of a full orthogonal group $O(N)$ or a symplectic
group $Sp(N)$. For these groups the
 defining or vector representations
have dimension $N$ (in the orthogonal case) and dimension $2N$
(in the symplectic case) respectively.

It is well-known that the isomorphism classes of simple representations
of $O(N)$ are labeled by Young diagrams with at most $N$ boxes in
the first two columns; simple representations of $Sp(N)$ are labeled
by Young diagrams with at most $N$ rows. We call such Young diagrams
permissible (for the respective group).

It is easy to describe
the decomposition of the tensor product of a simple representation
with the vector representation.
Let $X_\la$ be a simple object in $\Ca$ corresponding to the Young
diagram $\la$, and let $X=X_{[1]}$ be the object corresponding
to the vector representation (which is labeled by the Young diagram
with one box).
Then $X_\la\otimes X$ is the direct sum of simple representations
labeled by all permissible Young diagrams $\mu$ which are obtained from 
$\la$ by removing or adding a box from/to $\la$. While tensoring
with the vector representation would not per se describe the Grothendieck
semiring, it is all that we need for our purposes
together with the braiding (see Prop. \ref{Groth}).

In the following,  we denote by $[1^n]$ and  by $[n]$ the Young diagrams with
all its $n$ boxes in its first column and in its first row respectively.
The simple object
$X_{[1^n]}$  corresponds to the full antisymmetrization 
of the $n$-th tensor power $X^{\otimes n}$ of the vector representation
of the orthogonal group. In the representation category of symplectic
groups it would correspond to the unique simple subrepresentation in the
$n$-th antisymmetrization of the vector representation which has
not already appeared in the smaller tensor powers.
We obtain as a special case of the tensor product rule described
above
\begin{equation}\label{tens}
X_{[1^{m}]}\otimes X\ \cong\ X_{[1^{m+1}]}\oplus X_{[2,1^{m-1}]}
\oplus X_{[1^{m-1}]},\quad 1\leq m<N;
\end{equation}
if $m=N$, the right hand side above would be isomorphic to $ X_{[1^{N-1}]}$
in the orthogonal case, and to $ X_{[2,1^{N-1}]} \oplus X_{[1^{N-1}]}$
in the symplectic case.

\subsection{Fusion categories}\label{fusioncat}
There also exist braided tensor categories whose Grothendieck semirings
are quotients of the ones described in the last subsection. In these
cases, we can describe the labeling set for its simple objects
by also applying analogous restrictions to the rows of Young diagrams
as we had before for columns. We have the following three cases,
for fixed $N,M\in\N$:
\begin{enumerate}
\item orthogonal fusion category: the simple objects are labeled by
Young diagrams with $\leq N$ boxes in its first two columns and with
$\leq M$ boxes in its first two rows,

\item ortho-symplectic fusion category: the simple objects are labeled by
Young diagrams with $\leq N$ boxes in its first two columns and with
$\leq M$ boxes in its first row (i.e. with $\leq M$ columns),

\item symplectic fusion category: the simple objects are labeled by
Young diagrams with at most $N$ boxes in the first column and 
with at most $M$ boxes in the first row.
\end{enumerate}
Tensoring with the object labeled by the Young diagram with one
box (the analog of the `vector representation' in this context)
is as before, with now only those objects allowed at the right
hand side which satisfy the conditions for the labeling set of simple
objects in the corresponding fusion category. In particular, this
simple tensor product rule allows to compute the multiplicity of an
object $X_\la$ in $X^{\otimes n}=X_{[1]}^{\otimes n}$ by induction.

\subsection{Definition and examples}\label{examples}
In the rest of the paper, we have the following assumptions:
All categories are supposed to be
rigid, strictly monoidal, semisimple, braided $\C$-categories.
We say that such a category, say $\Ca$, is of $orthogonal$ or
\emph{symplectic type} if its Grothendieck semiring is the one of a
representation category of $O(N)$ (including $O(\infty)$) or
$Sp(N)$, or of one of the associated fusion categories, as described
in the last two subsections. Here are examples for such categories:

a) By definition, the representation categories 
$\rep\big(O(N)\big)$ and $\rep\big((Sp(N)\big)$ are tensor categories
of orthogonal resp. symplectic type, which have symmetric braiding.

b) It is well-known that the representation category of the
Drinfeld-Jimbo quantum group $U_q\g$ associated to the semisimple Lie
algebra $\g$ is semisimple  and that $\rep(U_q\g)$ has the same 
Grothendieck semiring as $\rep(\g)$, if $q$ is not a root of unity.  As
$\rep(\sp_N)$ is equivalent to $\rep\big(Sp(N)\big)$, $\rep(U_q\sp_N)$
is a braided tensor category of symplectic type. It is also possible
to construct braided tensor categories of orthogonal type as a
semidirect product of a subcategory of $\rep(U_q\so_N)$ with
$\rep(\Z/2)$. 

c) If $q$ is a root of unity, H.H. Andersen defined the category of
tilting modules of $U_q\g$. This category contains as a quotient
a semisimple category with only finitely many equivalence classes of
simple objects. These are examples of fusion categories. One can
construct the fusion categories of the last section from these
quotient categories in complete analogy to the construction sketched
in (b).

d) The existence of fusion categories was suggested by physicists
in conformal field theory. In particular, this implied the existence
of a highly nontrivial tensor product for representations of
affine Kac-Moody algebras resp. loop groups. A mathematically
rigorous definition was given by Kazhdan and Lusztig in
the Kac-Moody case (see \cite{kazhdan-lusztig1-2},
\cite{kazhdan-lusztig3},\cite{kazhdan-lusztig4})
and by Wassermann in \cite {wassermann-cft-3} for loop groups.
The equivalence between these categories and the ones defined by
Andersen was shown by Finkelberg \cite{finkelberg}.

e) It is also possible to construct orthogonal and symplectic categories
by topological methods
as quotients of the tangle category (see \cite{turaev-wenzl}).
This approach is closest to the set-up in this paper. It will be described
in more detail in Section \ref{existence}. A similar approach also
works for Lie type $A$ (see \cite{blanchet}).

\subsection{Low tensor powers} 
As $\ident$ is a subobject of $X^{\otimes 2}$, any simple subobject
of $X^{\otimes (n-2)}$ is also isomorphic to a simple subobject
of $X^{\otimes n}$. Hence we can write $X^{\otimes n}$
as a direct sum $X_{(n-2)}\oplus X_n$, where $X_{(n-2)}$
is a direct sum of simple objects each of which is isomorphic
to a subobject of $X^{\otimes (n-2)}$ and $X_n$ is a direct
sum of simple objects which are not isomorphic to any subobject of
$X^{\otimes (n-2)}$. By semisimplicity of $\Ca$,
we get from this the decomposition
\begin{equation}
\End\big(X^{\otimes n}\big)\cong \End\big(X_{(n-2)}\big)
\oplus \End\big(X_n\big).
\end{equation}

\begin{lemma}\label{linind1}
The set $\tilde\B =\{ 1, c, \tilde e\} \subset \End\big(X^{\otimes 2}\big)$ 
is linearly independent. In particular,  $c$ acts via different scalars
on $X_{[2]}$ and on $X_{[1^2]}$.

\end{lemma}

\begin{proof}
Assume  that $\tilde\B$ is not linearly independent. Then
we can assume $c=\alpha 1+\beta \tilde e$,
with $\alpha ,\beta\in F$, as otherwise the noninvertible $\e$ 
would be proportional to $1$. But then all the $c_i$'s just act
as scalars in End $(X_n)$. Let now $f$ resp. $\tilde f$ be the
projections onto the simple subobjects $X_{[1^2]}$ resp  $X_{[2]}$
of $X^{\otimes 2}$. Then we get, using $n=4$ and the braiding with
$c_{(2)}=c_{X^{\otimes 2},X^{\otimes 2}}=c_2c_1c_3c_2$
$$f\otimes \tilde f = f_1c_{(2)}\tilde f_1 c_{(2)}^{-1}
\quad \in \End\big(X_{(n-2)}\big),$$
where the last inclusion follows from the fact that $f_1\tilde f_1=0$
and that $c_{(2)}$ only acts as scalar in $\End(X_n)$, i.e. conjugation
by it induces the trivial automorphism in $\End(X_n)$.

As $\End\big(X_{[1^2]}\otimes X_{[2]}\big)
\cong f_1\tilde f_3 \, \End\big(X^{\otimes 4}\big)
\, f_1\tilde f_3\subset \End\big(X_{(2)}\big)$, we obtain that
$X_{[1^2]} \otimes X_{[2]}$ decomposes into a direct sum of simple 
modules which already appear in $X^{\otimes 2}$ (i.e. they are
isomorphic to $\ident , X_{[1^2]}$ or $X_{[2]}$); this contradicts
the tensor product rules for orthogonal and symplectic groups.
\end{proof}

\begin{lemma}\label{linind}
The space $\Hom\big(X,X^{\otimes 3}\big)$ has the basis
$\B=\left\{ i_2, c_1\circ i_2, \tilde e_1\circ i_2=i_1'\right\}$.
\end{lemma}

\begin{proof} This is a special case of Frobenius reciprocity:
the map $a\in \End\big( X^{\otimes 2}\big) \mapsto (a\otimes 1)\circ i_2$
has the inverse map $b\in\Hom\big( X,X^{\otimes 3}\big)\mapsto 
(1_2\otimes d)\circ (b\otimes 1)$. The claim now follows from Lemma
\ref{linind1}.
\end{proof}

\subsection{Matrix representations}\label{matrices} 
We define the quantity $d(X)$ by $d(X)= d\circ i$. Recall from the last
section that $d(X)=\al \dim(X)$ (see Lemma \ref{tracepreli}). 

\begin{lemma}\label{normalization}
There are scalars $r, q$ and a fourth root of unity
$\gamma$ such that 
\begin{enumerate}
\item the element $t=\gamma c$ has the eigenvalues $q$,
$-q^{-1}$ and $r^{-1}$ for the submodules
 $X_{[2]}$, $X_{[1^2]}$ and  $\1$ of $X^{\otimes 2}$ respectively,

\item if $q\neq q^{-1}$, then 

\begin{equation}\label{defdx}
d(X)=\frac{r-r^{-1}}{q-q^{-1}}+1=\frac{q^{-1}(r+q)(q-r^{-1})}{q-q^{-1}}.
\end{equation}

\item $\tr(t)=r/d(X)$.
\end{enumerate}
\end{lemma}

\begin{proof} 
It will be useful to compute matrix representations
of the elements $c_i$ and $\e_i$, $i=1,2$, acting on $\Hom(X,X^{\otimes 3})$
via concatenation. We will use the basis $\{ i_2, c_1\circ i_2, i_1\}$.
We claim that if the eigenvalues of $c$ are $\la_1,\la_2$ and $\la_3$,
then we obtain the matrix representations

\begin{equation}\label{matrices1}
 c_1\ \mapsto\ 
\begin{pmatrix}0&-\la_1\la_2&0\\
               1&\la_1+\la_2 &0\\
               0&\la_3(\la_1^{-1}+\la_2^{-1})&\la_3
\end{pmatrix}
\quad {\rm and }\quad 
 c_2\mapsto
\begin{pmatrix}\la_3&0&\la_1\la_2(\la_1+\la_2)\\
               0&0&-\la_1\la_2\\
               0&1&\la_1+\la_2
\end{pmatrix}.
\end{equation}
To see this observe that we have the obvious relations $c_j\circ i_j=
\la_3 i_j$ for $j=1,2$, and, from the braiding axiom,
$c_2\circ(c_1\circ i_2)=i_1$.
This determines two of the three columns of $c_j$, $j=1,2$.
Of the remaining column, two entries are computed using the fact that
the matrix must have determinant $\la_1\la_2\la_3$ and trace
$\la_1+\la_2+\la_3$. The remaining entries can be computed checking
the braid relation $c_1\circ c_2\circ c_1=c_2\circ c_1 \circ c_2$.
Moreover, using the braiding relation, we get $c_1\circ c_2\circ i_1=i_2$,
while the corresponding matrices, applied to $i_1$ would give 
$(\la_1\la_2)^2i_2$. Hence we also have 
$(\la_1\la_2)^2=1$, and we can assume $\la_1=\gamma^{-1} q$, $\la_2=-\gamma^{-1} q^{-1}$
and $\la_3=\gamma^{-1} r^{-1}$ for certain complex numbers $r$ and $q$ and for $\gamma$
a fourth root of unity.


Now it easily follows and the results of Lemma \ref{tracepreli},
\begin{equation}\label{matrices3}
\e_1\quad \mapsto\quad
\begin{pmatrix}0&0&0\\
               0&0&0\\
               \al&\al \la_3& \dim X
\end{pmatrix}
\quad {\rm and }\quad 
\e_2\quad\mapsto\quad
\begin{pmatrix}\dim X & \al\la_3^{-1}& \al\\
               0&0&0\\
               0&0&0
\end{pmatrix}.
\end{equation}
Comparing the $(3,2)$-matrix entries in the equality $\e_1c_1=\la_3\e_1$, we obtain
$$\big(\la_3(\la_1^{-1}+\la_2^{-1})\big)\big(\dim X\big)=\al\big(\la_3^2+\la_1\la_2-\la_3(\la_1+\la_2)\big).$$
If $\la_1+\la_2\neq 0$, this gives the formula for the dimension and
for $d(X)$ as stated,
after  substituting $r$ and $q$ into the eigenvalues as above.
It follows from this and Lemma \ref{tracepreli}, with
$\r=\al\la_3^{-1}$ that $\tr(t)=\tr(\gamma c)=r/d(X)$, as stated.

If $\la_1=-\la_2$, we  deduce from the last  equation  that $\la_3^2=-\la_1\la_2=
\pm 1=\la_2^2=\la_1^2$.
This implies that two of the three  eigenvalues of $c$
are identical and that the eigenvalues of $t$ are contained in the set $\{\pm 1\}$.
\end{proof}

\begin{lemma}\label{normallemma}
Let $t$ be as in Lemma \ref{normalization}. If $t$ has less than three
distinct eigenvalues, then necessarily its eigenvalues are $\pm 1$.
\end{lemma}

\begin{proof} We can rule out $\la_1=\la_2$ by Lemma \ref{linind1}.
Assume now that  $\la_1=\la_3$ or $\la_2=\la_3$, which would imply
$r=-q$ or $r=q^{-1}$ for  the eigenvalues of $t$.
If $\la_1+\la_2\neq 0$, we obtain  $\dim X=0$ from the computations  of the  previous 
lemma, which would contradict rigidity. If $\la_1+\la_2=0$, the claim follows 
from the last paragraph of the proof of the last lemma.
\end{proof}

\subsection{Relations}\label{relations} We can now summarize the results
of this section as follows: Let $e=i\circ d= \al \tilde e$.

\begin{proposition}\label{catrel}
\ 
\begin{enumerate}
\item Assume that $c$ has three distinct eigenvalues, and let $t=\gamma c$ be
as in Lemma \ref{normalization}. Then we can define the element 
$e\in\End\big(X^{\otimes 2}\big)$ also by $t-t^{-1}=(q-q^{-1})(1-e)$. We then 
have the relations
\begin{itemize}
\itemindent=.25in
\item[(R1)] $t_ie_i=r^{-1}e_i$, for $i=1,2,\ ...\ n-1$, and
\item[(R2)] $e_it_{i-1}^{\pm 1}e_i=r^{\pm 1}e_i$,  for $i=,2,\ ...\ n-1$.
\end{itemize}

\item If $c$ has fewer than three eigenvalues, then the representation of
the braid group $B_n$ given by the
morphisms $t_i$ factors through the symmetric group $S_n$. Moreover,
the elements $t_i$ and $e_i$, $i=1,2\ ...\ n-1$ generate a quotient
of Brauer's centralizer algebra.

\item We also have $\tr(t)=r/d(X)$ and $\tr(e)=1/d(X)$ and
$\tr\big((a\otimes 1)\chi_{n-1}\big)=\tr(a)\tr(\chi)$ for $\chi\in\{t,e\}$ in both cases;
here $\tr$ is the normalized trace on $\End \big(X^{\otimes}\big)$ and
$a\in \End\big(X^{\otimes n-1}\big)$.
\end{enumerate}
\end{proposition}

\begin{proof} By definition, $e$ is a multiple of the eigenprojection of $t$
for the  eigenvalue  $r^{-1}$. It can be seen e.g. from the explicit 
matrix representations, see \ref{matrices3}, that this multiple is $d(X)$.
The alternative formula for $e$ can now be checked easily, as well as
(R1). Part (c) follows from Lemma \ref{normalization} resp Lemma
\ref{tracepreli} for the values of $\tr(t)$ and $\tr(e)$, and
from Corollary \ref{condcor} for the Markov property. Using the
relation between $t^{\pm 1}$ and $e$ in part (a) of the statement,
one also computes $\tr(t^{-1})=r^{-1}/d(X)$.

By functoriality, it suffices to check Relation (R2) for $i=2$.
This follows from Lemma \ref{condlemma}(b) and (a), and from the values
of $\tr(t^{\pm 1})$ which have already been computed.
The proof for part (b) will be given in Section \ref{Brauer}.
\end{proof}

\section{$q$-Deformation of Brauer's centralizer algebra}

After having determined properties of braiding morphisms for braided
tensor categories of orthogonal or symplectic types, we now go
the  opposite way. We use the relations obtained in the last section
to  define  abstract algebras which turn out to be 
Brauer's centralizer algebras (see \cite{brauer}) and
a $q$-deformation of it which was discovered in connection with
Kauffman's link invariant (see \cite{birman-wenzl} and \cite{murakami-bmw}; here we 
follow the presentation in \cite{wenzl-bcd}, p 399/400).

\subsection{Hecke algebras}\label{Hecke} We first need a simpler class
of algebras. The Iwahori-Hecke algebra $\H_n(q^2)$ of type
$A_{n-1}$ is the algebra
defined over the field $F$
by generators $\Tt_i$, $i=1,2,\ ...\ n-1$, which satisfy
the braid relations 
and the quadratic relation
$\Tt_i^2=(q-q^{-1})\Tt_i+1$; here $q$ is a fixed element in
$F$. We have the following well-known theorem:

\begin{theorem}\label{Hecketheorem}
If $q^2$ is not a root of unity of order $\leq n$, then $\H_n(q^2)$
is isomorphic to the group algebra $F S_n$ of the symmetric group $S_n$.
\end{theorem}

One of the consequences of the last theorem is that the irreducible
representations of $\H_n(q^2)$ are labeled by Young diagrams with $n$
boxes if $\H_n(q^2)$ is semisimple. In that case, let 
$\tilde P_{[1^n]}$ 
be the central idempotent corresponding to the
one-dimensional representation $\Tt\mapsto -q^{-1}$.
Let $A\otimes 1$ denote the element in $\H_{n+1}$ obtained from
the element $A\in\H_n$ under the natural embedding of $\H_n$ into
$\H_{n+1}$. It is well-known that we have
\begin{equation}
\tilde P_{[1^n]}\otimes 1=\tilde P_{[1^{n+1}]} +
\tilde P_{[2,1^{n-1}]},
\end{equation}
where $\tilde P_{[2,1^{n-1}]}$ is an idempotent in the simple component
of $\H_{n+1}$ labeled by the Young diagram ${[2,1^{n-1}]}$.
Let $[n]_q=(q^n-q^{-n})/(q-q^{-1})$.

\begin{lemma}\label{Heckeidem} 
We have the following identities
in $\H_n$, for $m=1,2,\ ...\ n-1$:
\begin{enumerate}
\item $\tilde P_{[1^{m+1}]}=\frac{1}{[m+1]_q}
\left(q^m\tilde P_{[1^m]} -[m]_q \tilde P_{[1^m]}\Tt_m\tilde P_{[1^m]}\right)$
\item $\tilde P_{[1^m]}\Tt_m \tilde P_{[2,1^{m-2}]} \Tt_m\tilde P_{[1^m]}
=\frac{[m-1]_q[m+1]_q}{[m]_q^2} \left(\tilde P_{[1^m]}-\tilde P_{[1^{m+1}]}\right)
=\frac{[m-1]_q}{[m]_q} \tilde P_{[1^m]}\left(\Tt_m+q^{-1}1\right)\tilde P_{[1^m]} $.
\end{enumerate}
\end{lemma}

\begin{proof} These identities follow as special cases from properties of
path idempotents connected to Hoefsmit's orthogonal representations of Hecke
algebras (see e.g. \cite{wenzl-hecke}, Cor 2.3). They can also be proved
by induction on $m$ as follows: 

We can write $\Tt_i=(q+q^{-1})\tilde E_i- q^{-1}1$,
where $\tilde E_i$ is the eigenprojection for the eigenvalue $q$ of $\Tt_i$.
Then one  shows by induction on $m$, using $\tilde P_{[1^m]}\tilde E_{m-1}=0$
and $\tilde E_i\tilde E_{i-1}\tilde E_i=\tilde E_{i-1}\tilde E_i\tilde E_{i-1}
-\frac{[1]_q}{[2]_q^2}(\tilde E_i-\tilde E_{i-1})$
that
$$\tilde P_{[1^m]}\tilde E_{m}\tilde P_{[1^m]}\tilde E_{m}=
\frac{[m+1]_q}{[m]_q [2]_q}\tilde P_{[1^m]}\tilde E_{m}$$
and
$$\tilde P_{[1^{m+1}]}=
\tilde P_{[1^m]} - \frac{[m]_q [2]_q}{[m+1]_q}
\tilde P_{[1^m]}\tilde E_m\tilde P_{[1^m]}.$$
Claim (a) follows from the second equation.
Claim (b) follows by substituting $\tilde P_{[2,1^{m-1}]} =
\tilde P_{[1^{m+1}]}-\tilde P_{[1^{m}]}$, and then applying (a)
for $\tilde P_{[1^{m+1}]}$.
\end{proof}

\subsection{Definitions}\label{BMWdef} 
The algebra $\D_n(r,q)$, depending on two parameters $r$ and $q$,
is  given by generators $T_1,\ T_2
\ ...\ T_{n-1}$, which satisfy the braid relations 
and
\vskip .2cm

\begin{itemize}
\item[(R1)] $\quad E_iT_i=r^{- 1}E_i$,

\item[(R2)] $\quad E_iT_{i-1}^{\pm 1}E_i=r^{\pm 1}E_i,$

\vskip .2cm

\ni where $E_i$ is defined by the equation
\item[(D)] $ (q-q^{-1})(1-E_i)=T_i-T_i^{-1}.$
\end{itemize}
$Remarks:$ It is easy to read off from the defining relations the following
facts:

\begin{enumerate}
\item Let $\I_n$ be the ideal of $\D_n$ generated by $E_{n-1}$. Then 
$\D_n/\I_n\cong \H_n(q^2)$, with the isomorphism given
by $T_i\mapsto \tilde T_i$.

\item $\I_n\cong \D_{n-1}\otimes_{\D_{n-2}}\D_{n-1}$ as a $\D_{n-1} - \D_{n-1}$
bimodule, where the isomorphism is given by $D_1\otimes D_2 \mapsto
D_1E_{n-1}D_2$ for $D_1, D_2\in\D_{n-1}$.

\item If $\D_n$ is semisimple, $\D_n\cong \I_n\oplus \H_n$.

\item The $T_i$'s satisfy the cubic equation
$(T_i-r^{-1})(T_i+q^{-1})(T_i-q)=0$.
\end{enumerate}

\subsection{Structure of $q$-Brauer algebras}
The following Theorem
determines the structure of $\D_n(r,q)$ if it is semisimple
(see \cite{birman-wenzl}, Theorem 3.7 and \cite{wenzl-bcd}, Theorem 3.5 and Cor 5.6):

\begin{theorem}\label{BMWstruc}
\ 
\begin{enumerate}
\item The algebra $\D_n(r,q)$ is semisimple for generic values of $r$
and $q$ (see Theorem \ref{BMWtrace} for more specific information).
In this case, it has dimension $1\cdot 3\cdot 5\ ...\ (2n-1)$
and its simple components 
are labeled by the Young diagrams
with $n$, $n-2$, $n-4$, ..., 1 resp. 0
boxes.

\item The decomposition of a simple $\D_{n,\la}$ module $V_{n,\la}$ 
into simple $\D_{n-1}$ modules is given by
\begin{equation}
V_{n,\la}\cong \bigoplus_{\mu} V_{n-1,\mu},
\end{equation}
where the summation goes over all Young diagrams $\mu$ which can
be obtained by either taking away or, if $\la$ has less than $n$
boxes, by adding a box to $\la$. The labeling of simple components
is uniquely determined by this restriction rule and the convention
that the eigenprojection of $T_1$ corresponding to its eigenvalue
$q$ is labeled by the Young diagram $[2]$.

\item For diagrams $\la$ with $n$ boxes,  $V_{n,\la}$ 
becomes an $\H_n(q^2)$ module via the homomorphism of Remark (a)
in Section \ref{BMWdef}.

\item $\D_{n+1}= \spn \big\{ \ A\chi B \mid A,B\in \D_n,\ \chi\in \{ 1, T_n, E_n\} \big\}$.
\end{enumerate}
\end{theorem}

We leave it to the reader to check, using the inductive rule in Theorem
\ref{BMWstruc}, (b) (see also \cite{birman-wenzl}, Fig. 8) 
that $\D_1(r,q)\cong F$, $\D_2(r,q)\cong F^3$ and,
with $M_k(F)$ denoting the algebra of all $k\times k$ matrices,
$$\D_3(r,q)\cong M_3(F)\oplus F \oplus M_2(F)\oplus F.$$
It is an easy exercise to show (using relations (1)-(10) in
\cite{wenzl-bcd}, p.\ 400)  that the 3-dimensional simple component
contains a minimal left ideal spanned by $\{ E_2, T_1E_2, E_1E_2\}$,
and that the matrices which describe the action of the elements
$T_i$ and $E_i$, $i=1,2$ with respect to this basis
coincide with the ones in Eq \ref{matrices1} and 
\ref{matrices3}.

\subsection{Brauer algebras}\label{Brauer}
Brauer defined abstract finite dimensional algebras 
$\BD_n=\BD_n(x)$ (see \cite{brauer}) depending on a parameter $x$.
These abstract algebras are easiest
described by graphs. We will not give this description here (see \cite{brauer}).

The following description will suffice for our purposes: The algebras
$\BD_n$ can be defined via generators $T_i'$ and $E_i'$, $i=1,2,\ ...\ n-1$.
For $x=N>n$, we obtain a faithful surjective representation of $\BD_n(N)$ onto
$\End_{O(N)}\big(V^{\otimes n}\big)$ which maps $T_i'$ to the permutation of
the $i$-th and $(i+1)$-st factor, and it maps $E_i'$ to the element
$\e_i$ defined for this category as in Section \ref{seldual}; here
$V$ is the $N$-dimensional vector representation of $O(N)$.
Similarly, one obtains a surjective map $\BD_n(-2N)$
onto $\End_{Sp(N)}\big(V^{\otimes n}\big)$. 

The commutation relations between the
elements $T_i'$ and $E_i'$ are exactly the same ones as for the 
elements $T_i$ and $E_i$ in $\D_n$. In particular, the elements
$T_i'$, $E_i'$ commute with $T_j'$ and $E_j'$ whenever $|i-j|\geq 2$.
In fact, the relations for $x=N$ follow from the
ones in $\D_n(q^{N-1},q)$ in the limit for $q\rightarrow 1$.
(see e.g. \cite{birman-wenzl}, Section 5 or \cite{wenzl-bcd}, p 401 for details).

\begin{proof}[Conclusion of the  proof of Prop.\ref{catrel}:]
Evaluating the matrices in the proof of Lemma \ref{normalization} for
$r=q=1$, we obtain matrices for $t_i, e_i$, $i=1,2$ which only depend
on $d(X)$.  In particular $t_i^2=1$ for $i=1,2$. Moreover, if
$d(X)=N$, these matrices have to satisfy the same relations as the
corresponding elements in $\rep\big(O(N)\big)$. By functoriality, the
elements $t_i, e_i, t_{i+1}, e_{i+1}$ satisfy the same relations as
the elements $t_1,e_1,t_2,e_2$, and generators whose indices differ by
at least 2 commute. Hence the elements $t_i, e_i$ generate an algebra
isomorphic to a quotient of Brauer's centralizer algebra.
\end{proof}

\subsection{q-Dimensions}
We also need a general formula for $q$-dimensions of orthogonal
and symplectic groups.
Let $[y+n]_q=(rq^n-r^{-1}q^{-n})/(q-q^{-1})$.
Then we  define for each Young diagram $\la$
the rational function
\begin{equation}\label{defQ}
Q_{\lambda}(r,q)\ =
\ \prod_{(j,j)\in\lambda}\frac{(r-q^{-2\la_j+2j-1})(r^{-1}+q^{-2\la'_j+2j-1})}
{1-q^{-2h(j,j)}}
\quad \mathop{\prod_{(i,j)\in \lambda}}_{i\neq j}
\frac{[y+d(i,j)]_q}{[h(i,j)]_q};
\end{equation}
here $(i,j)$ denotes the box in the $i$-th row and $j$-th column
of $\la$,  $\la_i$ ($\la'_j$) is the number of boxes
in the $i$-th row ($j$-th column) of $\la$. Moreover, the quantity
$d(i,j)$ and the hook length $h(i,j)$  are defined by

\begin{equation}
d(i,j)\ = \begin{cases}
\lambda_i+\lambda_j-i-j+1 & \text{if $i\leq j$,}\\
-\lambda_i'-\lambda_j'+i+j-1 &\text{if $i> j$.}
\end{cases}
\end{equation}
and
\begin{equation}
h(i,j)=\lambda_i-i+\lambda_j' -j+1
\end{equation}
We will need these functions primarily for the special case of a Young
diagram $[1^m]$ whose only column contains exactly $m$ boxes. In this
case, we obtain
\begin{equation}\label{defQs}
Q_{[1^m]}(r,q)=\frac{(r-q^{-1})(r^{-1}+q^{1-2m})}{1-q^{-2m}}
\ \prod_{j=1}^{m-1} \frac{rq^{1-j}-r^{-1}q^{j-1}}{q^j-q^{-j}}.
\end{equation}
One checks similarly that for the Young diagram $[2,1^{m-2}]$ with
two boxes in the first row and one box in the next $m-2$ rows one obtain
\begin{equation}\label{defQss}
Q_{[2,1^{m-2}]}(r,q)=\frac{(r-q^{-3})(r^{-1}+q^{3-2m})}{1-q^{-2m}}
\frac{[y+1]_q [y+2-m]_q}{[1]_q [m-2]_q}
\prod_{j=1}^{m-3} \frac{rq^{1-j}-r^{-1}q^{j-1}}{q^j-q^{-j}}.
\end{equation}

The rational functions $Q_\la(r,q)$ have obvious analogs $\hat Q_\la(y)$
for the Brauer algebras. They are essentially defined by replacing
$q$-numbers in the definition of $Q_\la$ by ordinary numbers. More
precisely, we have
\begin{equation}\label{defQBr}
\hat Q_\la(y)=
\prod_{(i,j)\in \lambda}
\frac{y+d(i,j)}{h(i,j)},
\end{equation}

\subsection{Markov traces and semisimplicity}
The algebras $\D_n(r,q)$ carry 
an important trace functional  which we will describe 
in two different ways. The existence of the trace was originally
derived from the existence of Kauffman's link invariant
(see \cite{birman-wenzl},\cite{murakami-bmw}). The equivalent description
in the semisimple case follows from
\cite{wenzl-bcd}, Theorem 3.6 and Theorem 5.5. A more
algebraic existence proof can be given using the theory of 
quantum groups (see e.g. \cite{wenzl-bcd} and \cite{orellana-wenzl-spin-groups}, Lemma 3.1).

\begin{theorem}\label{BMWtrace}
\ 
\begin{itemize}
\item[(a1)] There exists a trace functional $\trd$ on $\D_n(r,q)$
which  is uniquely determined inductively by
$\trd(1)=1$, $\trd(T_i)=r/d(X)$, $\trd(E_i)=1/d(X)$
and by $\trd(A\chi B)=\trd(AB)\,\trd(\chi)$,
where $A,B\in D_{n-1}$ and $\chi\in\{ T_{n- 1}, E_{n- 1},1 \}$;
here $d(X)=(r-r^{-1})/(q-q^{-1})+1$ is defined as in Lemma 
\ref{normalization},(b).

\item[(a2)] If $\D_n$ is semisimple, the functional $\trd$ in (a) 
is completely determined by 
$\trd(p)=Q_\la/ d(X)^n$,
where $p$ is a minimal idempotent in $\D_{n,\la}$.

\item[(b)] Conversely, if $q^2$ is not a primitive $l$-th root of unity
for $1<l\leq n$, and if $Q_\la(r,q)\neq 0$ for all Young diagrams $\la$
with less than $n$ boxes, then the algebra $\D_n(r,q)$ is semisimple.

\item[(c)] If $r=q^{N-1}$,
$Q_\la(q^{N-1},q)$ is equal to the $q$-dimension
of the highest weight module $V_\la$ of $O(N)$. If $r=-q^{2N-1}$,
$(-1)^{|\la|}\,Q_\la(-q^{2N+1},q)$ is equal to the $q$-dimension
of the highest weight module $V_\la$  of $Sp(N)$, where $|\la|$ is the
number of boxes of $\la$. The 
$q$-dimension of $V_\la$ 
is defined to be the character of the element $q^{2\rho}$,
acting on $V_\la$, where $\rho$
is half the sum of the positive weights of the corresponding semisimple
Lie algebra.

\item[(d)] One can similarly define the Markov trace for the Brauer algebras
$\BD_n(d(X))$, where now the functions $Q_\la(r,q)$ are replaced by
the polynomials $\hat Q_\la(d(X))$ (with $r=q=1$).
\end{itemize}
\end{theorem}

\subsection{Quotients of $\D_n(r,q)$} \label{quotients}
It will be important to
compute the quotient of $\D_n(r,q)$ modulo the annihilator ideal
$\A_n$ of $\trd$, i.e. $\A_n=\{ A\in \D_n, \trd(AB)=0\ {\rm for\ all\ }
B\in\D_n\}$. In the following we assume $q^2$ to be a primitive $l$-th
root of unity (with $l=\infty$ covering the case $q^2=1$ or $q$ not
a root of unity).

\begin{itemize}
\item[(a1)] If $r=q^{N-1}$ or if $r=-q^{N-1}$ for $N$ odd,
with $q^2$ not a root of unity, then 
$\D_n/\A_n\cong \End_{O(N)}\big(V^{\otimes n}\big)$, where $V$ is the vector 
representation of the orthogonal group $O(N)$.
If $r=-q^{2N+1}$ with $q^2$ not a root if unity, then 
$\D_n/\A_n\cong \End_{Sp(N)}\big(V^{\otimes n}\big)$, where $V$ is the vector 
representation of the symplectic  group $Sp(N)$.
\item[(a2)] If $r$ is equal to $\pm$ a negative power of $q$ and $q^2$
is not a root of unity, then again $\D_n/\A_n$ is isomorphic to
$\End_G\big(V^{\otimes n}\big)$, with $V$ the vector representation of an
orthogonal or symplectic group $G$. The group can be determined 
from $(a1)$ after replacing $r$ by $-r^{-1}$.
The results listed in (a1) and (a2) are proved in \cite{wenzl-bcd}, Corollary 5.6.
\item[(b)] If  $q^2$ is a primitive $l$-th root of unity,
we can find positive integers $n,m<l$ such that $r=\pm q^n$ and $
r=\pm q^{-m}$ (where the signs may or may not match). Then we can
find restrictions for the number of boxes in the first (two) row(s) as
well as in the first (two) column(s), as it was described in parts (a1)
and (a2). Then again $\D_n /\A_n$ is isomorphic to $\End\big(X^{\otimes n}\big)$,
where now $X$ is the `vector representation' of the corresponding fusion
category, as described in Section \ref{fusioncat}.
See \cite{wenzl-bcd}, Theorem 6.4 for a somewhat more explicit description
and a proof.
\end{itemize}

\subsection{Reparametrization}\label{reparametrization}
It is easy to see that for the category
$\Ca$ generated by $X$, we have several different braiding structures.
Indeed, it is easy to check that replacing $c=c_{X,X}$ by $-c$, $c^{-1}$
or $-c^{-1}$ again gives a braiding structure. Moreover, we have made
a choice by labeling the object corresponding to the eigenvalue $q$
by the Young diagram $[2]$, and not by $[1^2]$. 
These observations
are reflected on the level of the algebras $\D_n(r,q)$ as follows:

\begin{enumerate}
\item There exist algebra isomorphisms $\D_n(r,q)\cong\D_n(-r,-q)
\cong \D_n(-r^{-1},-q^{-1})\cong$ \\
$\cong\D_n(r^{-1},q^{-1})$
given by
$$T_i\mapsto -T_i(-r,-q)\mapsto -T_i^{-1}(-r^{-1},-q^{-1})
\mapsto T_i^{-1}(r^{-1},q^{-1}),$$
where $T_i(r',q')$, $i=1,2,\ ...\ n-1$ are the generators of the algebra
$\D_n(r',q')$. These isomorphisms preserve the labeling of the simple 
components by Young diagrams.

\item 
There exists an isomorphism between $\D_n(r,q)$ and $\D_n(r,-q^{-1})$
by mapping $T_i(r,q)$ to $T_i(r,-q^{-1})$. This isomorphism maps
the simple component $\D_{n,\la}(r,q)$ to $\D_{n,\la'}(r,-q^{-1})$,
where $\la'$ is the Young diagram obtained from the Young
diagram $\la$ by interchanging rows with columns. 
By composing this isomorphism with the isomorphisms under (a), we obtain
additional isomorphisms which change the parametrization, e.g.
$\D_n(r,q)\cong \D_n(-r^{-1},q)$ is obtained by mapping $T_i(r,q)$
to $-T_i^{-1}(-r^{-1},q)$.

\item The isomorphisms in (a) and (b) preserve the Markov traces
(i.e. the pull-back of the Markov trace under one of these isomorphisms
gives the Markov trace of the original algebra). 

\item By uniqueness of the Markov trace, the isomorphisms in (a) and (b)
lead to identities for the functions $Q_\la(r,q)$ as follows:
$Q_\la(r,q)=Q_\la(-r,-q)=Q_\la(r^{-1},q^{-1})=
Q_{\la'}(r,-q^{-1})=Q_{\la'}(-r^{-1},q)$ etc.
\end{enumerate}
The statements above are easily proved (see also e.g \cite{wenzl-bcd}
Prop. 3.2(c)). It is also immediate that the isomorphisms above are
examples of {\it functorial isomorphisms} which are defined as
follows: Let $\bar \D_n(r,q)$ and $\bar \D_n(r',q')$ be quotients of
$\D_n(r,q)$ and $\D_n(r',q')$ respectively. We say that an isomorphism
$\Phi :\bar  \D_n(r,q)\to \bar \D_n(r',q')$ is $functorial$ if it maps
$\big\langle T_i(r,q)\big\rangle$ to $\big\langle
T_i'(r',q')\big\rangle$ for each $i$ with $1\leq i<n$; here $\langle
a\rangle $ is the subalgebra generated by an element $a$ of an algebra
$A$.

The following Lemma will result in another proof that the representation
category of $O(2)$ does not allow any deformations. We will denote by
$\bar \D_n(r,q)$ the quotient of $\D_n(r,q)$ with respect to the
annihilator ideal of $\trd$.

\begin{lemma}\label{O(2)} The algebras $\bar\D_n(q,q)$ and $\bar\D_n(q',q')$ are
functorially isomorphic for any $q,q'\in\C$ and any $n\in\N$.
\end{lemma}

\begin{proof} One checks easily that $Q_{[n]}(q,q)=2$ for $n>0$, that
$Q_{[0]}(q,q)=1=Q_{[1^2]}(q,q)$,and that $Q_\la(q,q)=0$ for all
other Young diagrams. One deduces from this that $\bar \D_n(q,q)\cong
\bar \D_n(q',q')$ as abstract algebras (see \cite{wenzl-bcd}, Cor 5.6(b3)).
In particular, $\bar D_3(q,q)$ is isomorphic to the direct sum of
a  full $3\times 3$ matrix algebra and a copy of $\C$. 
Let $p_i^{(\la)}$ be the eigenprojection of the element $t_i$ corresponding
to the object $X_\la$, with $\la\in\{ [0], [1^2],[2]\}$. Using the
basis $p_1^{(\la)}\circ i_2$ for $\Hom(X,X^{\otimes 3})$, one computes the
following matrices

\begin{equation}\label{matrices4}
 p_2^{([0])}\ \mapsto\ \frac{1}{4}
\ \begin{pmatrix}1&1&2\\
               1&1 &2\\
               1&1&2
\end{pmatrix},\ 
 p_2^{([1^2])}\ \mapsto\ \frac{1}{4}
\ \begin{pmatrix}1&1&-2\\
               1&1 &-2\\
               -1&-1&2
\end{pmatrix}
\text{, and }
 p_2^{([1^2])}\mapsto\ \frac{1}{2}
\ \begin{pmatrix}1&-1&0\\
               -1&1&0\\
               0&0&0
\end{pmatrix};
\end{equation}
this can be done fairly easily by using the dual basis
$\left\{ \frac{1}{\tr\big(p_1^{(\la)}\big)}\, d_2\circ p_1^{(\la)}\right\}$ with
$\la\in\{ [0], [1^2],[2]\}$ and the values for $Q_\la(q,q)$.
The crucial observation now is that these matrices do not depend
on $q$, and hence also the commutation relations between the
various  $p_1^{(\la)}$ and $p_2^{(\la)}$, modulo the annihilator
ideal of $\tr$. Hence we obtain $\bar\D_n(q,q)$ as the quotient
of an algebra whose defining relations are independent of $q$
with respect to the annihilator ideal of a trace functional
which does not depend on $q$ as well.
\end{proof}

\subsection{Inductive formulas for idempotents}\label{idformulas}
We will have to study the algebra $\D_n(r,q)$ for values of $r$ and $q$
for which it is not semisimple. This requires more explicit expressions
for certain central idempotents. These formulas are special cases
for inductive formulas of path idempotents, which have been studied
in \cite{ram-wenzl}. However, as we need somewhat more precise information,
including the nonsemisimple case,
we give a more or less self-contained derivation of the necessary results
here.

Let $A\in \D_m$. We shall denote 
by $A\otimes 1$ (or sometimes just by $A$, for brevity)
the image of $A$ under the
usual embedding of $\D_m$ into $\D_{m+1}$ which identifies the
generators of $\D_m$ with the first $m-1$ generators of $\D_{m+1}$.
Let $P_{[1^m]}$ denote the central idempotent belonging to
$D_{m,[1^m]}$ in the semisimple case. 
Using the restriction rule (2.1), we can write
\begin{equation}\label{(2.8)}
P_{[1^m]}\otimes 1  =P_{[1^{m+1}]} +\P2m + \Pmm1 ,
\end{equation}
where $\P2m$ is an idempotent in $\D_{m+1, [2,1^{m-1}]}$
and $\Pmm1$ is an idempotent in the simple component $D_{m+1, [1^{m-1}]}$.
By \cite{ram-wenzl}, (2.15), we have
\begin{equation}\label{(2.9)}
\Pmm1 = \frac{Q_{[1^{m-1}]}}{ Q_{[1^{m}]}}\ P_{[1^m]}  E_mP_{[1^m]}.
\end{equation}

\begin{lemma}\label{l2.5}
The idempotents $P_{[1^k]}$ are well-defined if $[m]_q\neq 0$ and
$r+q^{1-2m}\neq 0$ for $m=1,2,\ ...\ k$.
\end{lemma}

\begin{proof} 
The claim follows as soon as one has shown the following inductive
formula:
\begin{equation}\label{indform}
P_{[1^{m+1}]}=\frac{1}{[m+1]_q}\left(q^mP_{[1^m]}-[m]_qP_{[1^m]}T_mP_{[1^m]}
-\frac{[m]_q}{1+rq^{1-2m}}P_{[1^m]}E_mP_{[1^m]}\right).
\end{equation}
Observe that the algebra $\D_{m+1}$ is spanned by elements of the
form $A\chi B$, with $A,B\in \D_m$ and $\chi\in \{ 1, T_m,E_m\}$
(see Theorem \ref{BMWstruc},(d) or, e.g. \cite{wenzl-bcd}, Prop. 3.2). 
As $P_{[1^m]}A$ is a scalar multiple of
$P_{[1^m]}$ for any $A\in \D_m$, the subalgebra
$P_{[1^m]} \D_{m+1} P_{[1^m]}$ is spanned by the three elements
$P_{[1^m]} \chi P_{[1^m]}$, with  $\chi\in \{ 1, T_m,E_m\}$.
It follows from Lemma \ref{Heckeidem} and from $\D_n/\I_n\cong\H_n$
that we can write
$$
P_{[1^{m+1}]}= \frac{1}{[m+1]_q}\left(q^mP_{[1^m]}-[m]_q P_{[1^m]}T_mP_{[1^m]}
-\beta P_{[1^m]}E_m P_{[1^m]}\right)
$$
for some suitable scalar $\beta$. To compute the scalar,
we evaluate each side of the equation above under $\trd$. Using
the Markov property, we obtain
$$
\frac{Q_{[1^{m+1}]}}{d(X)^{m+1}}=
\frac{1}{[m+1]_q}
\left(\frac{q^mQ_{[1^m]}}{d(X)^{m}}-\frac{[m]_q r Q_{[1^m]}}{d(X)^{m+1}}
-\beta\frac{ Q_{[1^m]}  }{d(X)^{m+1}}\right).
$$
Using the explicit formula for $Q_{[1^m]}$ (see Eq \ref{defQs}),
one can easily solve for $\beta$.
\end{proof}

\begin{lemma}\label{l2.6}
Assume that $r=q^{m-1}$, with $m>0$ and that $q^2$ is a primitive
$l$-th root of unity, $l>m+1$ or $l=\infty$. Then
$P_{[1^{m+1}]}$ is well-defined and
$P_{[1^{m+1}]}\otimes 1$  is a  central minimal idempotent
in $D_{m+2}$ modulo the ideal $\J$ generated by $\P2m$.
\end{lemma}

\begin{proof} It is easy to check that the expressions for $P_{[1^k]}$
in Lemma \ref{l2.5} are well-defined for our choice of parameters
if $k\leq m+1$; this also implies that $\P2m$ is well-defined.
As $\P2m$ is a linear combination of $P_{[1^m]}\chi P_{[1^m]}$,
with $\chi\in \{ 1,T_m, E_m\}$, it follows from the relations
that $E_{m+1}\P2m E_{m+1}$ is a scalar multiple of $E_{m+1}P_{[1^m]}$.
The scalar can be computed to be equal to $Q_{[2,1^{m-1}]}/Q_{[1^m]}$
by using the Markov property of $\trd$. 
Using this, one easily shows that  
$P_{[1^{m+1}]} E_{m+1}P_{[1^{m+1}]}\in \J$.
As  $\D_{m+2}/\I_{m+2}\cong\H_{m+2}$, it follows from Lemma \ref{Heckeidem}
and from $\D_n/\I_n\cong \H_n$ that 
\begin{equation}\label{ideal}
P_{[1^{m+1}]} T_{m+1}\P2m T_{m+1} P_{[1^{m+1}]} = 
\frac{[m]_q}{[m+1]_q} P_{[1^{m+1}]}(T_{m+1}+q^{-1}) P_{[1^{m+1}]}
+ \gamma P_{[1^{m+1}]} E_{m+1}P_{[1^{m+1}]},
\end{equation}
where $\gamma$ is some scalar. This implies that also 
$ P_{[1^{m+1}]}(T_{m+1}+q^{-1}) P_{[1^{m+1}]}$ is in $\J$.
This shows that $P_{[1^{m+1}]}\otimes 1\equiv
P_{[1^{m+2}]}$ mod $\J$, if the latter is well-defined.

If $q^2$ is a primitive $(m+2)$-nd root of unity, we choose as
spanning set for the subalgebra $P_{[1^{m+1}]}\D_{m+2}P_{[1^{m+1}]}$ the
elements $P_{[1^{m+1}]}$, $P_{[1^{m+1}]}(T_{m+1}+q^{-1})P_{[1^{m+1}]}$
and $P_{[1^{m+1}]}E_{m+1}P_{[1^{m+1}]}$ and show as before that the 
last two elements are in $\J$.
\end{proof}

\begin{lemma}\label{nilpotent}
Let $q^2$ be a primitive $l$-th root of unity and assume 
$Q_{[1^k]}(r,q)\neq 0$ for $1\leq k\leq l$. Then there exists
a nilpotent element $N_l\in \D_l(r,q)$ such that $\trd(N_l)\neq 0$.
\end{lemma}

\begin{proof} We see from Lemma \ref{l2.5} that the elements $P_{[1^k]}$
are well-defined for $k<l$, and that also $N_l={[l]_q}P_{[1^l]}$
is well-defined. It follows that $N_l^2=[l]_qN_l=0$ for our choice of $q$.
Moreover, we have
$$\trd(N_l)=[l]_q \, \frac{Q_{[ 1^{l}]}}{ d(X)^{l}}.$$
It is easy to see from Eq \ref{defQs}
that $Q_{[1^l]}$ has a pole of order 1 for
our choice of $q$, which cancels with the zero of $[l]_q$ 
in the formula above.
Hence $\trd(N_l)\neq 0$ also for $q^2$ a primitive $l$-th root of unity.
\end{proof}

\begin{corollary}\label{nilcor}
The algebra $\D_l /\A_l$ is not semisimple if $q^2$ is a primitive $l$-th
root of unity and $Q_{[l]}(r,q)\neq 0$ or $Q_{[1^l]}\neq 0$.
\end{corollary}

\begin{proof} If $Q_{[1^l]}(r,q)\neq 0$, we can find an element $N_l\in
\D_l(r,q)$ which has nonzero trace (hence also nonzero in the quotient
mod $\A_l$) but it is nilpotent. This is not possible
in a semisimple algebra. The case with $Q_{[l]}(r,q)$ goes similarly,
using one of the isomorphisms in Section \ref{reparametrization}.

The quotient $\D_n/\A_n$ is semisimple for all $n\in\N$
if and only if
$Q_{[l]}(r,q)=0$ and $Q_{[1^l]}(r,q)=0$; this condition is vacuous 
for $l=\infty$. The `only if' part follows from Corollary \ref{nilcor}.
The `if' part follows from below where we list all the other
cases for the parameters $r$ and $q$.
\end{proof}

\section{Identifying $\End\big(X^{\otimes n}\big)$}

We have seen in the last two sections that there exists a homomorphism
$\Phi$ from the 
algebra $\D_n(r,q)$ or $\BD_n(d(X))$ into $\End\big(X^{\otimes n}\big)$
given by $T_i\mapsto t_i$ and $E_i\mapsto e_i$.
The purpose of this section is to show that this map is surjective.

\subsection{Preliminaries}
We say that two idempotents $e$ and $f$ in an algebra
$\M$ are (von Neumann) equivalent, $e \sim f$, if there exist
elements $u$ and $v$ in $\M$ such that $e=uv$ and $f=vu$.
An idempotent $e\in\M$ is called minimal if there exists for any
$a\in \M$ a scalar $\gamma(a)$ such that $eae=\gamma(a)e$.
The \emph{multiplicity} $\mult_\M(e)$ of an  idempotent 
$e\in \M$ is the maximum
number $m$ of idempotents $e_i\in \M$ such that $e_ie_j=0$ for $i\neq j$
and $e_i \sim e_j$.

Recall that in a semisimple category we can associate to a subobject $X$ of
an object $Y$ (i.e. a monomorphism from $X$ into $Y$) an idempotent $p_X$ in
$\End(Y)$ (see e.g. Lemma \ref{quasiinverse}). We then define the 
multiplicity of the subobject $X$ in $Y$ to be equal to the multiplicity of
$p_X$ in $\End(Y)$.

\begin{lemma}\label{3.2}
Let $\Ca$ be a semisimple category, let $Y\in \Ob(\Ca)$ and let
$e,f\in$ $\End(Y)$ be two idempotents.
\begin{enumerate}
\item The idempotents $e$ and $f$ are equivalent
iff $\im(e)$ and $\im(f)$ are isomorphic subobjects of $Y$.

\item Let $X$ be a subobject of $\im(e)$. Then the multiplicity
of $X$ in $Y$ is $\geq \mult_{\End(Y)}(e)$.
\end{enumerate}
\end{lemma}

\begin{proof} Follows straightforward from the definitions.
\end{proof}

\subsection{} Let now $\Ca$ be a (fusion) category of orthogonal or
symplectic type, and let $N$ be the maximum of numbers $k$ for which
we have a simple object in $\Ca$ labeled by a Young diagram of the
form $[1^k]$.

\begin{lemma}\label{3.3} Let $1\leq m<N$ and assume that $P_{[1^{m}]}$ and 
$P_{[1^{m+1}]}$ exist in $\D_{m+1}(r,q)$.
Let $p_{[1^{k}]}=\Phi\left(P_{[1^{k}]}\right)$ for $k=1,2,\ ...,\ m+1$.
\begin{enumerate}
\item If $\im\left(p_{[1^{m}]}\right)=X_{[1^{m}]}$, then $X_{[1^{m+1}]}$ is a subobject
of $\im\left(p_{[1^{m+1}]}\right)$.

\item If moreover $\Phi\left(\P2m\right)\neq 0$, then  $\im\left(p_{[1^{m+1}]}\right)=X_{[1^{m+1}]}$.
\end{enumerate}
\end{lemma}

\begin{proof} If $m=1$, the statements are true by definition. Assume now
$m>1$. By induction assumption,
$X_{[1^m]}$ is a subobject of $\Phi\left(P_{[1^{m}]}\right)$; hence
 $X_{[1^{m+1}]}$ is a subobject of $\Phi\left(P_{[1^{m}]}\otimes 1\right)$.
As $X_{[1^{m+1}]}$ has multiplicity 1 in $X^{\otimes (m+1)}$,
the claim in (a) follows from  Eq \ref{(2.8)} and Lemma \ref{3.2}.(b).

For part (b), observe that $(p_{[1^{m}]}\otimes 1)
\End\big(X^{\otimes m+1}\big)(p_{[1^{m}]}\otimes 1)$ has dimension 3 by
Eq \ref{tens}. On the other hand, 
$\trd (\Pmm1)=\dim X_{[1^{m-1}]}/(\dim X)^{m+1}\neq 0$. Using this,
our assumption on  $\Phi\left(\P2m\right)$ and part (a), it follows that
the three idempotents
on the right hand side of Eq \ref{(2.8)}
have nonzero image under $\Phi$.
Hence the claim follows from Eq \ref{tens}.
\end{proof}

\subsection{Restrictions for parameters} Let $\Ca$ and $N$ be as
in the previous subsection. Recall that we can choose a fourth root of
unity $\gamma$ such that the 
eigenvalues of $\gamma c$ are $q$, $-q^{-1}$ and $r^{-1}$ for suitable
values $q$ and $r$.

\begin{lemma}\label{3.4}
Assume that $q^2$ is a primitive $l$-th root of unity,
$l\in \N\cup\{ \infty \}$ and let $m\in\N$ be such 
that $Q_{[1^{m+1}]}(r,q)=0$ and $Q_{[1^k]}(r,q)\neq 0$ for
$1\leq k\leq m$.
Then $m<l$ and $m\leq N$.
\end{lemma}

\begin{proof} 
Assume $l\leq m$. By Lemma \ref{nilpotent}, there exists a nilpotent
element $N_{l}\in\D_{l}$ with $\trd(N_{l})\neq 0$.  Then also
$\Phi(N_{l})$ is nilpotent and $\tr\big(\Phi(N_{l})\big)=\trd(N_l)\neq
0$, a contradiction to $\End\big(X^{\otimes l}\big)$ being semisimple.

Now assume that $m>N$. Then it follows from Lemma \ref{l2.5} that
$P_{[1^{N+1}]}$, $P_{[2,1^{N-1}]}$ and $P_{[1^{N-1}]}^{(N+1)}$
are well-defined. By our assumptions, they also
have nonzero trace. From this we could conclude that
$\Phi\big((P_{[1^N]}\otimes 1)(\D_{N+1})(P_{[1^N]}\otimes 1)\big)$
would have dimension $\geq 3$. This contradicts the fact
that dim $\End\big(X_{[1^N]}\otimes X\big)=1$ in the orthogonal case
and dim $\End\big(X_{[1^N]}\otimes X\big)=2$ in the symplectic case (see
the remark below Eq \ref{tens}).
\end{proof}

\begin{lemma}\label{3.5}
\ 
\begin{enumerate}
\item If $\Ca$ has the Grothendieck semiring of an orthogonal group
$O(N)$ or of one of its associated fusion categories,
then $r=q^{N-1}$ or, if $N$ is odd, $r=-q^{N-1}$.

\item If  $\Ca$ has the Grothendieck semiring of a symplectic group
$Sp(N)$ or of one of its associated fusion categories,
then $r=-q^{2N+1}$.
\end{enumerate}
\end{lemma}

\begin{proof} Let $m$ be as in Lemma \ref{3.4}. Assume $m<N$.
If $\Phi\big(P_{[2,1^{m-1}]}\big)\neq 0$, then 
$\im \Phi\big(P_{[1^{m+1}]}\big)=X_{[1^{m+1}]}$ by Lemma \ref{3.3} and 
$\dim X_{[1^{m+1}]} = Q_{[1^{m+1}]}(r,q)=0$, which contradicts
rigidity, Lemma \ref{rigidity}.

If $\Phi\big(P_{[2,1^{m-1}]}\big)=0$, then $X_{[2,1^{m-1}]}$ has to be
a subobject of $W=$ $\im\big(\Phi(P_{[1^{m+1}]})\big)$ by Eq
\ref{(2.8)} ($\im\big(\Phi(P_{[1^{m-1}]}^{(m+1)})\big)$ is
isomorphic to $X_{[1^{m-1}]}$); in particular, $W\cong X_{[1^{m+1}]}
\oplus X_{[2,1^{m-1}]}$ is not a simple object.  Moreover,
$\big(P_{[2,1^{m-1}]}\otimes 1\big) \subset \ker \Phi$ and
$\Phi\big(P_{[1^{m+1}]}\otimes 1\big)$ is a central and minimal
idempotent in $\Phi(\D_{m+2})$ by Lemma \ref{l2.6}.  By the braiding
axioms, we can identify $c_{W,X}$ with an element in
$\Phi\big(P_{[1^{m+1}]}\otimes 1\big)\D_{m+2}\big(P_{[1^{m+1}]}
\otimes 1\big) \cong \C$. Hence $c_{W,X}$ is a scalar multiple of
$1_{W\otimes X}$.  As
$$c_{W, X^{\otimes n}}=
\big(1_{ X^{\otimes n-1}}\otimes c_{W,X}\big) \circ 
\big(c_{W\otimes X^{\otimes n-1}}\otimes 1_X\big),$$
it follows that $c_{W, X^{\otimes n}}$ is a multiple of the
identity for all $n\in\N$. As $W$ is a subobject of $X^{\otimes m+1}$,
we also get that $c_{W,W}$ is a multiple of $1_{W\otimes W}$.
But then conjugation by $c_{W,W}$ would not permute the factors of
$p_{[1^{m+1}]}\otimes p_{[2,1^{m-1}]}\subset \End(W)\otimes \End(W)$,
with $p_{[1^{m+1}]}$, $p_{[2,1^{m-1}]}$ the projections onto the
submodules of $W$, contradicting the braiding property.
This, together with Lemma \ref{3.4} forces $m=N$.

Using the formulas \ref{defQs}, one checks that $m=N$ implies
$r=q^{N-1}$, $r=-q^{N-1}$ if $N$ is odd, or $r=-q^{2N+1}$.
In case (a) we can rule out $r=-q^{2N+1}$, as in this case also
$Q_{[2,1^{N-1}]}\big(-q^{2N+1},q\big)\neq 0$. In case (b), we can rule out
the other cases for $r$ by observing that this would imply
$\dim X_{[2,1^{N-1}]} = Q_{[2,1^{N-1}]}\big(-q^{2N+1},q\big)=0$, which would 
contradict rigidity.
\end{proof}

\subsection{} We can now prove the main result of this section

\begin{theorem}\label{surjectivity}
Let $\Ca$ be a  tensor category of orthogonal or symplectic
type. Then the map $\Phi: \D_n(r,q)\to \End\big(X^{\otimes n}\big)$
induced by $T_i\mapsto t_i$ and $E_i\mapsto e_i$
is a well-defined, surjective algebra homomorphism, 
with the kernel being the annihilator
ideal $\A_n$ of the trace $\trd$. 
\end{theorem}

\begin{proof} We have seen in the proof of Lemma \ref{3.5} that a restriction
on the number of antisymmetrizations forces $r$ to be equal to $\pm$
a positive power of $q$. Similarly, it follows from the results in
Section \ref{reparametrization}
that a restriction on the number of symmetrizations forces 
$r$ to be equal to $\pm$ a negative power of $q$. In particular,
if we have restrictions of both the numbers  of symmetrizations and
antisymmetrizations, the two resulting
equalities force $q$ to be a root of unity. 
It now follows from  Section \ref{quotients} that the
quotient of $\D_n(r,q)$ modulo the annihilator ideal
of the categorical trace coincides with $\End\big(X^{\otimes n}\big)$.
\end{proof}
As an application of this theorem, we can now show that the description
of orthogonal and symplectic categories in Section \ref{fusioncat}
was sufficient.

\begin{proposition}\label{Groth} The Grothendieck semiring of a
category of orthogonal
or symplectic type is already uniquely determined by the labeling
set of its simple objects and the tensor
product rules involving the vector representation, see Sections
\ref{combdat} and \ref{fusioncat}.
\end{proposition}

\begin{proof} Observe that in all our paper we have only used the
tensor product rules involving the vector representation to prove
the last  theorem. By that theorem, any simple object $X_\la$
corresponds to an idempotent $p_\la$ in a quotient $\bar\D_n(r,q)$
of $\D_n(r,q)$ for some $n\in\N$.  With the simple object $X_\mu$
corresponding to an idempotent $p_\mu\in \bar\D_m(r,q)$,
the multiplicity of $X_\nu$ in $ X_\la\otimes X_\mu$ is now equal
to the multiplicity of the idempotent $p_\la\otimes p_\mu$ 
in the simple component of $\bar\D_{n+m}(r,q)$ labeled by $\nu$.

It only remains to show that the multiplicity of this idempotent
does not depend on the values of the parameters $r$ and $q$
for the various cases (see Sections
\ref{combdat}, \ref{fusioncat} and \ref{quotients}).  A 
proof probably most suited to our setting goes as follows:
A set of minimal idempotents for the algebra $\bar\D_n$ was defined in
\cite{ram-wenzl}, Cor. 2.5 (see also Section \ref{idformulas}).
Strictly speaking, this was only done there for
the generic case when $\D_n$ is semisimple, but the proof goes through
exactly the same way for $\bar\D_n$.  More precisely, inductive expressions
were given in terms of the generators with coefficients 
being rational functions
in $r$ and $q$ whose singularities are contained in the set of zeros of the 
dimension functions $Q_\la(r,q)$ for our given category. 
Moreover,  explicit matrix
representations were determined for the generators of the algebra
$\bar\D_n(r,q)$ whose matrix entries again are rational functions with 
singularities as before, see \cite{leduc-ram}, Theorem 6.15.

If $\bar\D_n(r,q)\cong\bar\D_n(r',q')$ for all $n$
and we are not in  the case of a fusion category,
we can  find a path $(r(t),q(t))$, $0\leq t\leq 1$ from $(r,q)$ to $(r',q')$
for which $\bar\D_n(r,q)\cong\bar\D_n(r(t),q(t))$ for  $0\leq t\leq 1$,
avoiding any possible pole for the matrix representing the idempotent
$p_\la\otimes p_\mu$.
By continuity, the rank of this idempotent hence
must be constant in each irreducible representation of
$\bar\D_n$ if we vary the parameters $r$ and $q$ along our
chosen path. For showing the claim
in the case of fusion categories, we can find a Galois isomorphism
which maps $(r,q)$ to $(r',q')$ (after possibly using some of the
reparametrizations mentioned in Section \ref{reparametrization}).
This again leaves the rank invariant.
\end{proof}

$Remarks$: 
The argument in the last proposition works as well in other cases where
the braiding elements of a generating object $X$ of a braided category
generate $\End\big( X^{\otimes n}\big)$ for all $n\in\N$. In particular,
it can be used for Lie type $A$ and the associated fusion
rings (see \cite{kazhdan-wenzl}).

\section{Classification of categories of orthogonal or symplectic types}

Let $\Ca$ be a tensor category of orthogonal or symplectic type, and let
$r$ and $q$ be the parameters deduced from the eigenvalues of the braiding
morphism $c_{X,X}$, see Lemma \ref{normalization}. We will show that these
parameters will essentially uniquely determine $\Ca$, up to a few 
special cases.

\subsection{Special Cases}\label{except} Let us first rule out a few
cases for which the following general discussion will not apply:
Observe that these include all the possible values of 
the parameters $r$ and $q$ for which $Q_{[2,1]}(r,q)=0$ (see Eq \ref{defQss}).
\begin{enumerate}
\item It is not possible that $q$ is a root of unity and $r$ is not a root
of unity. In this case we would obtain a nilpotent element $a$ in
$\End\big(X^{\otimes n}\big)$ for some $n\in\N$ with nonzero categorical trace,
which would contradict semisimplicity of $\End\big(X^{\otimes n}\big)$
(see Lemma \ref{nilpotent} and its corollary).

\item It is not possible that  $r=q^{-1}$ or $r=-q$; this would imply $d(X)=0$,
contradicting rigidity of $\Ca$.

\item It is not possible that  $r=\pm 1$ and $q\neq \pm 1$. In this case
 $Q_{[1^2]}(1,q)=0=Q_{[2]}(1,q)$, which would contradict rigidity

\item If $r=q$ or $r=-q^{-1}$ (the $O(2)$-case), we obtain a 
unique description
of $\End\big(X^{\otimes n}\big)$ independent of any parameters $r$ and $q$
(see Lemma \ref{O(2)}). Hence the diagonal $\dd$ in the $O(2)$ case
does not depend on $q$, and there exist exactly two monoidal algebras
in this case by Theorem \ref{monclass}.

\item If $r=q^{-3}$ or $r=-q^3$ (the $Sp(1)$-case), $Q_{[1^2]}$ resp
$Q_{[2]}$ is equal to 0. Hence in this case we can only obtain a rigid
category for which the braiding morphism for the object $X$ has only
two distinct eigenvalues. Such categories have been
classified  in \cite{kazhdan-wenzl} and, for this special case,
already before in \cite{frohlich-kerler}.
\end{enumerate}

\subsection{Existence}\label{existence}
We have already seen examples of orthogonal or
symplectic tensor categories in Section \ref{examples}. The most natural
construction in our context uses the tangle category (see 
\cite{joyal-street-braided}, \cite{yetter}
\cite{kassel}, \cite{turaev}). For more details about this construction
see \cite{turaev-wenzl} and, for the classical case, \cite{deligne}.

An $(n,m)$-tangle is a collection of $(n+m)/2$ ribbons and an arbitrary
number of annuli in $\R^2\times [0,1]$; moreover, $n$ ends of the ribbons
will be in the intervals $[i-\ep, i+\ep]\times \{ 0\}\times \{ 0\}$, 
$i=1,2,\ \dots, n$, and $m$ ends of the ribbons will be in the intervals
$[j-\ep, j+\ep]\times \{ 0\}\times \{ 1\}$, $j=1,2,\ \dots, m$. 
The concatenation $t_1\circ t_2$ of an
$(m,k)$-tangle $t_1$ with an $(n,k)$-tangle $t_2$ is given by putting 
$t_1$ on top of $t_2$
and rescaling the $z$-coordinate.

We want to use tangles to construct monoidal algebras. In order to
get finite dimensional morphism spaces, we need some relations between
various tangles. These are the Kauffman skein relations (see 
\cite{kauffman}, or also e.g. \cite{wenzl-bcd}). To do so consider
the following (0,2) and (2,0) tangles

\begin{picture}(260,60)(0,0)
\put(100,35){\qbezier(0,0)(15,-20)(30,2)}
\put(135,35){}
\put(110,13){$\iota$}
\put(200,35){\qbezier(0,0)(15,20)(30,-2)}
\put(235,35){}
\put(210,13){$\pi$}
\put(150,0){Figure 1}
\end{picture}

Here one should think of the ribbon obtained by thickening the lines
parallel to the drawing plane.
Then $\iota\circ \pi$ is a (2,2) tangle. Further (2,2) tangles  are given
by 1 (two parallel vertical ribbons) and $\sigma^{\pm 1}$
(two crossing ribbons, where the $\pm 1$ exponent corresponds to the two
possible ways of crossing them). We have two possible ways of defining
quotients, via the Kauffman skein relations:
\begin{equation}\label{kauff1}
\sigma -\sigma^{-1}=(q-q^{-1})(1-\iota\circ \pi)\quad
{\rm and}\quad
\sigma\circ\iota = r^{-1}\iota
\end{equation}
or
\begin{equation}\label{kauff2}
\sigma+\sigma^{-1}=\sqrt{-1}(q-q^{-1})(1+\iota\circ \pi)\quad
{\rm and}\quad
\sigma\circ\iota = \sqrt{-1}r^{-1}\iota.
\end{equation}
One can show that the $\C$-span of (0,0) tangles modulo these relations
is isomorphic to $\C$. 
Using this and the morphisms $\iota$ and $\pi$ similar as
the morphisms $i$ and $d'$ in Section 2, one defines a trace $\tr$ on
the set of $(n,n)$-tangles
(see e.g. chapters 2 and 3 in \cite{turaev-wenzl}).
A $\C$-linear combination $a$
of $(n,m)$-tangles is called \emph{negligible} if $\tr(ab)=0$ for any
$(m,n)$-tangle $b$.
Let $\bar\Ta(n,m)_{\pm}$ be the quotient of the 
$\C$-span of all $(n,m)$-tangles
modulo the negligible linear combinations of $(n,m)$ tangles with
respect to the trace defined by
relations \ref{kauff1} (for +) or \ref{kauff2} (for $-$).
Then it can be checked that, for chosen sign, the collection 
$(\bar\Ta(n,m)_{\pm})_{n,m\in\N}$ is a monoidal algebra
of type $2$ with  $\bar\Ta(n,n)_{\pm}\cong \D_n(r,q)/\A_n$, whenever the latter
is semisimple for all $n\in \N$.
From these monoidal algebras, one can construct semisimple categories
using the results of  Section \ref{catreconstruction}.
This has already been shown before in \cite{turaev-wenzl}, Theorem 8.6.
So we have
\begin{proposition}\label{exist} 
There exist categories of orthogonal or symplectic types
as quotient categories of the tangle category 
modulo relations \ref{kauff1} or \ref{kauff2} for all values $r,q$
for which $\bar\Ta(n,n)_{\pm}\cong \D_n(r,q)/\A_n$ is semisimple 
for all $n\in\N$. These cases
are all listed in Section \ref{quotients}.
\end{proposition} 
If $q=\pm 1$ (resp $q=\pm i$ in case of relation \ref{kauff2}), 
one only obtains interesting categories if
also $r=\pm 1$ (resp $r=\pm  i$); otherwise $d(X)$ would not be well-defined.
Moreover, one needs to add to these relations the additional relation 
$\pi\circ\iota = d(X)$,
with $d(X)\in \C$. In this case, it  is often more convenient to consider the
resulting structure as a category of graphs (see the work of Brauer
\cite{brauer} and Deligne \cite{deligne}. Then one obtains monoidal 
algebras and tensor
categories as in the previous proposition (see \cite{deligne}).
Moreover, using the polynomials \ref{defQBr}, one obtains 
(see \cite{wenzl-brauer}, Cor 3.3 and Cor. 3.5)

\begin{proposition}\label{exist2}
There exist orthogonal and symplectic categories
obtained as quotient 
categories of the tangle (or graph) category if $q\in\{\pm 1\}$ or
$q\in\{\pm i\}$, with the additional relation $\pi\circ\iota = d(X)$.
The resulting category $\Ca$ has the 
Grothendieck semiring of $\rep\big(O(\infty)\big)$
if $d(X)$ is not an integer. If $d(X)$ is an integer, $\Ca$ has the
Grothendieck semiring of $\rep\big(O(N)\big)$ if $d(X)=N$ or,
if $N$ is odd, $d(X)=2-N$
and it has the Grothendieck semiring of $Sp(N)$ if $d(X)=-2N$.
\end{proposition}

\subsection{Uniqueness} Let $\Ca$, $\tilde\Ca$ be categories of orthogonal 
or symplectic type with isomorphic Grothendieck semirings,
and let $r,q$ resp $\tilde r, \tilde q$ be the corresponding parameters
as determined in Lemma \ref{normalization}. We can rule out the special cases
considered in Section \ref{except}; in particular we can assume that
$Q_{[2,1]}$ is not zero for these  parameters.
Let $X$ and $\tilde X$
be objects corresponding to the (analogue of the) vector representation
in $\Ca$ and  $\tilde\Ca$ respectively.
Also, recall that as $X$ generates $\Ca$, its braiding structure
is uniquely determined by $c_{X,X}$.

\begin{theorem}\label{uniqueness} Let the notation be as above,
and assume that  $q\not\in\{ \pm 1\}$.
Then $\tilde \Ca$ is equivalent to
$\Ca$ as monoidal categories if and only if the eigenvalues of 
$ c_{\tilde X,\tilde X}$ can be obtained
from the ones of $c_{X,X}$
by changing the braiding and/or the labeling as described  in parts
(a) and (b) in Section \ref{reparametrization}.

If $q=\pm 1$, then categories $\Ca$ and $\tilde\Ca$ constructed as in 
Prop. \ref{exist2} are equivalent if and only if $d(X)= d(\tilde X)$
for the additional parameters  $d(X)$ and $ d(\tilde X)$, and
$c_{X,X}$ and $ c_{\tilde X,\tilde X}$ have the same eigenvalues.
\end{theorem}

\begin{proof}
Let $p^{(\la)}$ be the eigenprojection of $t$ for 
$X_\la$ is a subobject of $X^{\otimes 2}$. It is a well-known result for
Hecke algebras of type $A$, that the nonzero eigenvalue of 
$p^{(\la)}_1p^{(\la)}_2p^{(\la)}_1$ in the summand labeled by $[2,1]$ is equal to
$(q+q^{-1})^{-2}$ (see e.g.
\cite{wenzl-hecke}, p. 361).
Hence if $\tilde\Ca$ is equivalent to $\Ca$, we obtain
$(q+q^{-1})^{-2}=(\tilde q+\tilde q^{-1})^{-2}$, which entails 
$\tilde q\in \{ \pm q^{\pm 1}\}$. Hence, after changing the braiding 
structure in $\tilde \Ca$ by replacing $c_{X,X}$ by its negative
and/or inverse, if necessary, we can assume $\tilde q=q$.
It also follows that the quantities $d(X)^2$ and $Q_{[1^2]}$ must be the same 
for $\Ca$ and $\tilde\Ca$. Hence we obtain
\begin{equation}\label{d(X)}
\frac{\tilde r-\tilde r^{-1}}{q-q^{-1}}+1=\pm \left(\frac{r-r^{-1}}{q-q^{-1}}+1\right).
\end{equation}
If we have a plus sign on the right hand side, it follows that $\tilde r
\in \{r, -r^{-1}\}$, as claimed. To exclude the minus sign, one uses
$Q_{[1^2]}(r,q)=Q_{[1^2]}(\tilde r,q)$ 
(see Eq. \ref{defQss}) as follows: After substituting
the factor $(\tilde r-\tilde r^{-1})/(q-q^{-1})$, using Eq \ref{d(X)},
one obtains a second equation in which the only powers of $\tilde r$
are $\tilde r$ and $\tilde r^{-1}$. 
Solving this linear system in unknowns $\tilde r$ and $\tilde r^{-1}$,
it would follow that $\tilde r$
is a rational function of $r$ and $q$. However, this is not possible
for the solution of the quadratic equation \ref{d(X)} (in $\tilde r$);
it is easy to find integer values for $r$ and $q$ for which $\tilde r$
is not rational. This finishes the proof of one direction.

On the other hand, assume we have
orthogonal or symplectic categories $\Ca$ and $\tilde\Ca$ with isomorphic
Grothendieck semirings,
with the parameters $(r,q)$ and $(\tilde r, \tilde q)$ related as in
the statement. Hence, after suitable relabeling and change of braiding
structure, if necessary, we can assume that the braiding elements $c_{X,X}$
and $ c_{\tilde X,\tilde X}$ have the same eigenvalues, 
for the same components. By Theorem \ref{surjectivity},
this means that both $\End\big(X^{\otimes n}\big)$
and  $\End\big(\tilde X^{\otimes n}\big)$ are isomorphic
to $\bar\D_n(r,q)=\D_n(r,q)/\aaa_n$ for all $n\in \N$. Moreover, under this
isomorphism, the tensor operations in $\Ca$ and $\tilde\Ca$ correspond
to the usual embeddings of $\bar\D_n(r,q)\otimes \bar\D_m(r,q)$ into
$\bar\D_{n+m}(r,q)$. Hence we obtain an equivalence of the diagonal
monoidal algebra generated by $X$ and $\tilde X$.
By Theorem \ref{Theta-surjective} and its corollary, this equivalence
extends to the  monoidal algebras  generated by  $X$ and $\tilde X$.
 But then also $\Ca\cong\tilde\Ca$ by
Theorem \ref{reconstruction}. This completes the proof of the theorem 
if $q\neq \pm 1$.

As the quantity $d(X)$ is independent of the choice of $\iota$ and $\pi$,
equivalent categories of symplectic or orthogonal type  must have
the same value for $d(X)$. On the  other hand, if $q=\pm 1$ for two
categories $\Ca$ and $\tilde \Ca$ of orthogonal or symplectic type
for which also $d(X)=d(\tilde X)$, their diagonal  monoidal  algebras
are given by the Brauer algebras with parameter $d(X)= d(\tilde X)$,
hence are isomorphic.
As before, their two possible extensions can be told apart by the eigenvalues
of the braiding morphism $c_{X,X}$, by Corollary \ref{monclass}.

\end{proof}

\subsection{Main Theorem}  Let $\Ca$ be a tensor category of orthogonal 
or symplectic type, and let $X$ be the object corresponding to the
vector representation.

\begin{theorem}
\ 
\begin{enumerate}
\item The category $\Ca$ is  completely determined, up to the symmetries
mentioned in  Theorem \ref{uniqueness}, by the eigenvalues of the 
braiding
morphism $c_{X,X}$, which can be assumed to be of the form $q,-q^{-1}$ and 
$r^{-1}$ or of the form  $iq,-iq^{-1}$ and $ir^{-1}$, and if $q\in \{\pm 1,
\pm i\}$, by
the quantity $d(X)=\pi\circ \iota$.

\item The category $\Ca$ is a fusion category if and only if $q$ is a root of
unity and $r=\pm q^n$ for some $n\in\Z$ (see Section \ref{fusioncat});
it is of $O(N)$ or $Sp(N)$ type if and only if $r=\pm q^n$
with $n$ as in Section \ref{combdat} and $q$ not a root of unity
or if $q=\pm 1$ and
$d(X)$ is an integer,  and it is of type $O(\infty)$
if and only if $r$ is not $\pm$ a power of $q$ and $q$ is not a root of
unity. Moreover, such categories exist for all possible
values of $r$ and $q$, subject to these conditions, which have not already
been excluded in  Section \ref{except}.
\end{enumerate}
\end{theorem}

\begin{proof} Part (a) follows from Theorem \ref{uniqueness}.
Part (b) now follows from Theorem \ref{surjectivity}
and the results listed in Section \ref{quotients}; the existence
part follows from Propositions \ref{exist} and \ref{exist2}.

\end{proof}
Let $\E$ be a fundamental domain for the $\Z/2\times \Z/2$ action on 
$\C\backslash \{0\}$ given by $q\to q^{-1}$ and $q\to -q$. 

\begin{corollary}\label{catclass}
Braided tensor categories whose Grothendieck semirings are isomorphic to the one
of $\rep\big(Sp(N)\big)$ are in 1-1 correspondence with pairs $(q,\epsilon)$ where $q$
is a complex number  in $\E$ not  equal to a root  of unity except $\pm 1$ and 
$\epsilon\in \{\pm 1\}$. The same holds if $Sp(N)$ is replaced by  $O(N)$
with $N$ even. For odd $N$, we have two families  of braided tensor categories
each of  which  is labelled by pairs $(q,\epsilon)$ as above, which correspond
to the cases with $r=q^{N-1}$ and $r=q^{1-N}$.
\end{corollary}

\begin{proof}
By Theorem  \ref{surjectivity} it suffices to determine all  pairs of parameters
$(r,q)$ for which $\D_n(r,q)/\A_n\cong \End\big( X^{\otimes  n}\big)$ for all
$n\in\N$. Using the symmetries in Section \ref{reparametrization} and the results
in \cite{wenzl-bcd},  Theorem 6.4 
(see Section \ref{quotients}), one shows first
that we can assume the parameters $(r,q)$ to be of the form $(q^m,q)$,  with $q\in\E$.
Again  using \cite{wenzl-bcd}, Theorem 6.4 
(and \cite{wenzl-brauer}, Cor. 3.5 for $q=1$),
one can read  off  which exponent belongs to which group.
\end{proof}

$Remark:$ The categories whose Grothendieck  semirings  are isomorphic to the ones
of  a symplectic or an even-dimensional orthogonal group as well as  one of the
two families in the odd-dimensional orthogonal case are closely related  to
the corresponding  Drinfeld-Jimbo quantum groups. The second family of categories
in the odd-dimensional case seems to be different. For instance, it is not possible
to obtain positive dimensions for all objects,  for any choice of parameters,
even after changing the quantity $\alpha$ (see
Lemma \ref{tracepreli}) for the dimension function.
\bigskip

\nocite{andersen}
\nocite{birman-braids-links}
\nocite{birman-wenzl}
\nocite{bruguieres}
\nocite{finkelberg}
\nocite{murakami-bmw}
\nocite{jimbo}
\nocite{kac}
\nocite{kassel}
\nocite{kazhdan-lusztig1-2}
\nocite{kazhdan-lusztig3}
\nocite{kazhdan-lusztig4}
\nocite{kazhdan-wenzl}
\nocite{shnider}
\nocite{andersen-paradowski}
\nocite{blanchet}
\nocite{drinfeld-cocommutative}
\nocite{joyal-street-braided}
\nocite{kirillov}
\nocite{leduc-ram}
\nocite{lusztig-quantum-groups}
\nocite{moore-seiberg}
\nocite{morton-wassermann}
\nocite{turaev-wenzl}
\nocite{reshetikhin-quantum-groups}
\nocite{vogel}
\nocite{wassermann-cft-3}
\nocite{wenzl-c*}
\nocite{weyl-classical-groups}
\nocite{witten-qft-jones-poly}
\nocite{andersen-tensor-categories}
\nocite{bakalov-kirillov}
\nocite{frohlich-kerler}
\nocite{verlinde}
\nocite{wenzl-braidsurvey}

\bibliographystyle{abbrv}
\bibliography{bibliography}

\def\cprime{$'$}
\begin{thebibliography}{10}

\bibitem{andersen}
H.~H. Andersen.
\newblock Tensor products of quantized tilting modules.
\newblock {\em Comm. Math. Phys.}, 149(1):149--159, 1992.

\bibitem{andersen-tensor-categories}
H.~H. Andersen.
\newblock Quantum groups, invariants of {$3$}-manifolds and semisimple tensor
  categories.
\newblock In {\em Quantum deformations of algebras and their representations
  (Ramat-Gan, 1991/1992; Rehovot, 1991/1992)}, volume~7 of {\em Israel Math.
  Conf. Proc.}, pages 1--12. Bar-Ilan Univ., Ramat Gan, 1993.

\bibitem{andersen-paradowski}
H.~H. Andersen and J.~Paradowski.
\newblock Fusion categories arising from semisimple {L}ie algebras.
\newblock {\em Comm. Math. Phys.}, 169(3):563--588, 1995.

\bibitem{bakalov-kirillov}
B.~Bakalov and A.~Kirillov, Jr.
\newblock {\em Lectures on tensor categories and modular functors}, volume~21
  of {\em University Lecture Series}.
\newblock American Mathematical Society, Providence, RI, 2001.

\bibitem{birman-braids-links}
J.~S. Birman.
\newblock {\em Braids, links, and mapping class groups}.
\newblock Princeton University Press, Princeton, N.J., 1974.
\newblock Annals of Mathematics Studies, No. 82.

\bibitem{birman-wenzl}
J.~S. Birman and H.~Wenzl.
\newblock Braids, link polynomials and a new algebra.
\newblock {\em Trans. Amer. Math. Soc.}, 313(1):249--273, 1989.

\bibitem{blanchet}
C.~Blanchet.
\newblock Hecke algebras, modular categories and {$3$}-manifolds quantum
  invariants.
\newblock {\em Topology}, 39(1):193--223, 2000.

\bibitem{brauer}
R.~Brauer.
\newblock On algebras which are connected with the semisimple continuous
  groups.
\newblock {\em Ann. of Math.}, 63:854--872, 1937.

\bibitem{bruguieres}
A.~Brugui\`eres.
\newblock Reconstruction des certaines cat\'egories monoidales d'apr\`es
  {K}azhdan-{W}enzl.
\newblock Notes for seminar talks. University of Paris.

\bibitem{deligne}
P.~Deligne.
\newblock La s\'erie exceptionnelle de groupes de {L}ie.
\newblock {\em C. R. Acad. Sci. Paris S\'er. I Math.}, 322(4):321--326, 1996.

\bibitem{drinfeld-cocommutative}
V.~G. Drinfel{\cprime}d.
\newblock Almost cocommutative {H}opf algebras.
\newblock {\em Algebra i Analiz}, 1(2):30--46, 1989.

\bibitem{finkelberg}
M.~Finkelberg.
\newblock An equivalence of fusion categories.
\newblock {\em Geom. Funct. Anal.}, 6(2):249--267, 1996.

\bibitem{freyd}
P.~Freyd.
\newblock {\em Abelian categories. {A}n introduction to the theory of
  functors}.
\newblock Harper \& Row Publishers, New York, 1964.

\bibitem{frohlich-kerler}
J.~Fr{\"o}hlich and T.~Kerler.
\newblock {\em Quantum groups, quantum categories and quantum field theory},
  volume 1542 of {\em Lecture Notes in Mathematics}.
\newblock Springer-Verlag, Berlin, 1993.

\bibitem{jimbo}
M.~Jimbo.
\newblock A $q$-difference analogue of ${U}({\mathfrak g})$ and the
  {Y}ang-{B}axter equation.
\newblock {\em Lett. Math. Phys.}, 10(1):63--69, 1985.

\bibitem{joyal-street-braided}
A.~Joyal and R.~Street.
\newblock Braided tensor categories.
\newblock {\em Adv. Math.}, 102(1):20--78, 1993.

\bibitem{kac}
V.~G. Kac.
\newblock {\em Infinite-dimensional {L}ie algebras}.
\newblock Cambridge University Press, Cambridge, third edition, 1990.

\bibitem{kassel}
C.~Kassel.
\newblock {\em Quantum groups}, volume 155 of {\em Graduate Texts in
  Mathematics}.
\newblock Springer-Verlag, New York, 1995.

\bibitem{kauffman}
L.~H. Kauffman.
\newblock An invariant of regular isotopy.
\newblock {\em Trans. Amer. Math. Soc.}, 318(2):417--471, 1990.

\bibitem{kazhdan-lusztig1-2}
D.~Kazhdan and G.~Lusztig.
\newblock Tensor structures arising from affine {L}ie algebras. {I}, {I}{I}.
\newblock {\em J. Amer. Math. Soc.}, 6(4):905--947, 1993.

\bibitem{kazhdan-lusztig3}
D.~Kazhdan and G.~Lusztig.
\newblock Tensor structures arising from affine {L}ie algebras. {I}{I}{I}.
\newblock {\em J. Amer. Math. Soc.}, 7(2):335--381, 1994.

\bibitem{kazhdan-lusztig4}
D.~Kazhdan and G.~Lusztig.
\newblock Tensor structures arising from affine {L}ie algebras. {I}{V}.
\newblock {\em J. Amer. Math. Soc.}, 7(2):383--453, 1994.

\bibitem{kazhdan-wenzl}
D.~Kazhdan and H.~Wenzl.
\newblock Reconstructing monoidal categories.
\newblock In {\em I. M. Gel\cprime fand Seminar}, volume~16 of {\em Adv. Soviet
  Math.}, pages 111--136. Amer. Math. Soc., Providence, RI, 1993.

\bibitem{kirillov}
A.~A. Kirillov, Jr.
\newblock On an inner product in modular tensor categories.
\newblock {\em J. Amer. Math. Soc.}, 9(4):1135--1169, 1996.

\bibitem{leduc-ram}
R.~Leduc and A.~Ram.
\newblock A ribbon {H}opf algebra approach to the irreducible representations
  of centralizer algebras: the {B}rauer, {B}irman-{W}enzl, and type {A}
  {I}wahori-{H}ecke algebras.
\newblock {\em Adv. Math.}, 125(1):1--94, 1997.

\bibitem{lusztig-quantum-groups}
G.~Lusztig.
\newblock {\em Introduction to quantum groups}, volume 110 of {\em Progress in
  Mathematics}.
\newblock Birkh\"auser Boston Inc., Boston, MA, 1993.

\bibitem{maclane-cat}
S.~Mac~Lane.
\newblock {\em Categories for the working mathematician}.
\newblock Springer-Verlag, New York, second edition, 1998.

\bibitem{moore-seiberg}
G.~Moore and N.~Seiberg.
\newblock Classical and quantum conformal field theory.
\newblock {\em Comm. Math. Phys.}, 123(2):177--254, 1989.

\bibitem{morton-wassermann}
H.~R. Morton and A.~Wassermann.
\newblock A basis for the {B}irman-{W}enzl algebra.
\newblock 1999/2000.
\newblock Preprint, University of Liverpool.

\bibitem{murakami-bmw}
J.~Murakami.
\newblock The {K}auffman polynomial of links and representation theory.
\newblock {\em Osaka J. Math.}, 24(4):745--758, 1987.

\bibitem{orellana-wenzl-spin-groups}
R.~C. Orellana and H.~G. Wenzl.
\newblock {$q$}-centralizer algebras for spin groups.
\newblock {\em J. Algebra}, 253(2):237--275, 2002.

\bibitem{ram-wenzl}
A.~Ram and H.~Wenzl.
\newblock Matrix units for centralizer algebras.
\newblock {\em J. Algebra}, 145(2):378--395, 1992.

\bibitem{reshetikhin-quantum-groups}
N.~Y. Reshetikhin.
\newblock Quantized universal enveloping algebras, the yang-baxter equation and
  invariants of links.
\newblock 1988.
\newblock Preprint LOMI E-4-87.

\bibitem{shnider}
S.~Shnider.
\newblock Deformation cohomology for bialgebras and quasi-bialgebras.
\newblock In {\em Deformation theory and quantum groups with applications to
  mathematical physics (Amherst, MA, 1990)}, volume 134 of {\em Contemp.
  Math.}, pages 259--296. Amer. Math. Soc., Providence, RI, 1992.

\bibitem{Tuba-Wenzl-Braid}
I.~Tuba and H.~Wenzl.
\newblock Representations of the braid group ${B}\sb 3$ and of ${\rm
  {s}{l}}(2,{\bf {z}})$.
\newblock {\em Pacific J. Math.}, 197(2):491--510, 2001.

\bibitem{turaev-wenzl}
V.~Turaev and H.~Wenzl.
\newblock Semisimple and modular categories from link invariants.
\newblock {\em Math. Ann.}, 309(3):411--461, 1997.

\bibitem{turaev}
V.~G. Turaev.
\newblock {\em Quantum invariants of knots and 3-manifolds}, volume~18 of {\em
  de Gruyter Studies in Mathematics}.
\newblock Walter de Gruyter \& Co., Berlin, 1994.

\bibitem{verlinde}
E.~Verlinde.
\newblock Fusion rules and modular transformations in {$2$}{D} conformal field
  theory.
\newblock {\em Nuclear Phys. B}, 300(3):360--376, 1988.

\bibitem{vogel}
P.~Vogel.
\newblock Algebraic structures on modules of diagrams.
\newblock 1997.
\newblock Preprint posted at \texttt{www.math.jussieu.fr/\\ \~{}vogel}.

\bibitem{wassermann-cft-3}
A.~Wassermann.
\newblock Operator algebras and conformal field theory. {III}. {F}usion of
  positive energy representations of {${\rm LSU}(N)$} using bounded operators.
\newblock {\em Invent. Math.}, 133(3):467--538, 1998.

\bibitem{wenzl-hecke}
H.~Wenzl.
\newblock Hecke algebras of type {$A\sb n$} and subfactors.
\newblock {\em Invent. Math.}, 92(2):349--383, 1988.

\bibitem{wenzl-brauer}
H.~Wenzl.
\newblock On the structure of {B}rauer's centralizer algebras.
\newblock {\em Ann. of Math. (2)}, 128(1):173--193, 1988.

\bibitem{wenzl-bcd}
H.~Wenzl.
\newblock Quantum groups and subfactors of type {$B$}, {$C$}, and {$D$}.
\newblock {\em Comm. Math. Phys.}, 133(2):383--432, 1990.

\bibitem{wenzl-c*}
H.~Wenzl.
\newblock {$C^*$} tensor categories from quantum groups.
\newblock {\em J. Amer. Math. Soc.}, 11(2):261--282, 1998.

\bibitem{wenzl-braidsurvey}
H.~Wenzl.
\newblock Tensor categories and braid representations.
\newblock In {\em Quantum groups and Lie theory (Durham, 1999)}, volume 290 of
  {\em London Math. Soc. Lecture Note Ser.}, pages 216--234. Cambridge Univ.
  Press, Cambridge, 2001.

\bibitem{weyl-classical-groups}
H.~Weyl.
\newblock {\em The classical groups}.
\newblock Princeton Landmarks in Mathematics. Princeton University Press,
  Princeton, NJ, 1997.
\newblock Their invariants and representations, Fifteenth printing, Princeton
  Paperbacks.

\bibitem{witten-qft-jones-poly}
E.~Witten.
\newblock Quantum field theory and the {J}ones polynomial.
\newblock {\em Comm. Math. Phys.}, 121(3):351--399, 1989.

\bibitem{yetter}
D.~N. Yetter.
\newblock Markov algebras.
\newblock In {\em Braids (Santa Cruz, CA, 1986)}, volume~78 of {\em Contemp.
  Math.}, pages 705--730. Amer. Math. Soc., Providence, RI, 1988.

\end{thebibliography}
\end{document}